\documentclass[10pt]{amsart}
\usepackage{fancyvrb,amssymb,epic,eepic,verbatim,graphicx}
\input psfig.sty

\newtheorem{theorem}{Theorem}[subsection]

\newtheorem{corollary}[theorem]{Corollary}

\newtheorem{proposition}[theorem]{Proposition}
\newtheorem{lemma}[theorem]{Lemma}
\newtheorem{lem}[theorem]{}
\theoremstyle{definition}
\newtheorem{definition}[theorem]{Definition}
\theoremstyle{remark}
\newtheorem{remark}[theorem]{Remark}
\newtheorem{example}[theorem]{Example}
\newcommand{\blem}{\begin{lem} \rm}
\newcommand{\elem}{\end{lem}}
\newcommand\M{\mathcal{M}}

\renewcommand\M{\mathcal{M}}

\renewcommand{\L}{\mathcal{L}}

\newcommand{\R}{\mathbb{R}}

\newcommand{\RR}{\mathcal{R}}
\newcommand{\FF}{\mathcal{F}}
\newcommand{\C}{\mathbb{C}}

\newcommand{\Y}{\mathcal{Y}}
\newcommand{\Z}{\mathbb{Z}}

\newcommand{\ppthot}{\frac{\partial}{\partial \theta_{13}}}
\newcommand{\ppthtf}{\frac{\partial}{\partial \theta_{24}}}

\renewcommand{\P}{\mathbb{P}}

\newcommand\lie[1]{\mathfrak{#1}}

\newcommand{\m}{\lie{m}}

\renewcommand{\sl}{\lie{sl}}
\newcommand{\on}{\operatorname}

\newcommand{\Hom}{ \on{Hom}}

\newcommand{\Vol}{  \on{Vol}}
\newcommand{\diag}{  \on{diag}}

\newcommand{\ssm}{\kern-.5ex \smallsetminus \kern-.5ex}

\newcommand\dirac{/\kern-1.2ex\partial} 
\newcommand\qu{/\kern-.7ex/} 
\newcommand\lqu{\backslash \kern-.7ex \backslash} 

\newcommand\dr{r_+ \kern-.7ex - \kern-.7ex r_-}
 
\newcommand{\labell}\label

\renewcommand{\d}{{\mbox{d}}}
\newcommand{\ol}{\overline}

\newcommand\eps{\epsilon}

\newcommand{\f}{\frac}

\newcommand{\hh}{{\f{1}{2}}}
\newcommand{\qq}{{\f{1}{4}}}
\newcommand{\thh}{{\f{3}{2}}}

\newcommand\Tr{\on{Tr}}

\newcommand\Ai{\on{Ai}}
\newcommand\Bi{\on{Bi}}

\newcommand\bdefn{\begin{definition}}
\newcommand\edefn{\end{definition}}
\newcommand\bea{\begin{eqnarray*}}
\newcommand\eea{\end{eqnarray*}}
\newcommand\bcv{\left[ \begin{array}{r} }
\newcommand\ecv{\end{array} \right] }

\newcommand\bma{\left[ \begin{array} }
\newcommand\ema{\end{array} \right]}
\newcommand\bsj{\left\{ \begin{array}{rrr} }
\newcommand\esj{\end{array} \right\}}

\newcommand\pplab{\frac{\partial}{\partial l_{ab}}}
\newcommand\ben{\begin{enumerate}}
\newcommand\een{\end{enumerate}}
\newcommand\bex{\begin{example}}
\newcommand\eex{\end{example}}
\newcommand\id{\, \setlength{\unitlength}{0.00010in}
{
\begin{picture}(24,2439)(0,500)
\path(12,2412)(12,12)
\end{picture}
}
\,}
\renewcommand{\Y}[3]{\setlength{\unitlength}{0.00010in}
{
\begin{picture}(2424,3039)(0,500)
\path(12,2412)(1212,1212)(2412,2412)
\path(1212,1212)(1212,12)
{\tiny \put(12,2412){\makebox(0,0)[lb]{$#2$}}
\put(1812,2412){\makebox(0,0)[lb]{$#3$}}
\put(1212,12){\makebox(0,0)[lb]{$#1$}}}
\end{picture}
}
}
\newcommand\YYL[5]{\setlength{\unitlength}{0.00020in}
{
\begin{picture}(3024,2439)(0,500)
\path(12,2412)(1212,1212)(2412,2412)
\path(1212,2412)(612,1812)
\path(1212,1212)(1212,12)
{\tiny \put(112,2312){\makebox(0,0)[lb]{$#3$}}
\put(1212,2212){\makebox(0,0)[lb]{$#4$}}
\put(2412,2212){\makebox(0,0)[lb]{$#5$}}
\put(912,1512){\makebox(0,0)[lb]{$#2$}}
\put(1212,12){\makebox(0,0)[lb]{$#1$}}}
\end{picture}}}
\newcommand{\YYR}[5]	{\setlength{\unitlength}{0.00020in}
{
\begin{picture}(3024,2439)(0,500)
\path(12,2412)(1212,1212)(2412,2412)
\path(1212,1212)(1212,12)
\path(1212,2412)(1812,1812)
{\tiny \put(112,2312){\makebox(0,0)[lb]{$#3$}}
\put(1312,2312){\makebox(0,0)[lb]{$#4$}}
\put(2412,2212){\makebox(0,0)[lb]{$#5$}}
\put(1412,1212){\makebox(0,0)[lb]{$#2$}}
\put(1212,12){\makebox(0,0)[lb]{$#1$}}}
\end{picture}}}

\newcommand\lE{l_{cd}}

\newcommand\sx{*\kern-.5ex_X}

 \newcommand\Area{\on{Area}}

\def\mathunderaccent#1{\let\theaccent#1\mathpalette\putaccentunder}
\def\putaccentunder#1#2{\oalign{$#1#2$\crcr\hidewidth \vbox
to.2ex{\hbox{$#1\theaccent{}$}\vss}\hidewidth}}

\begin{document}


\title{$6j$ symbols for $U_q(\sl_2)$ and non-Euclidean tetrahedra}

\author{Yuka U. Taylor}

\author{Christopher T. Woodward} \thanks{Partially supported by NSF
grant DMS0093647 } 

\address{Department of Mathematics,
The George Washington University,
Old Main, Room 102,
  1922 F St. NW,
Washington, DC 20052}

\address{Mathematics-Hill Center,
Rutgers University, 110 Frelinghuysen Road, Piscataway, NJ 08854-8019,
U.S.A.}  \email{ctw@math.rutgers.edu}

\begin{abstract} We relate the semiclassical asymptotics of
the $6j$ symbols for the quantized enveloping algebra
$U_q(\lie{sl}_2)$ at $q$ a root of unity (resp. $q$ real positive) to
the geometry of spherical (resp. hyperbolic) tetrahedra.
\end{abstract}

\maketitle

\vskip -2in

\section{Introduction}

Let $r > 2$ be an integer, $q = \exp( \pi i/r)$ and $U_q(\lie{sl}_2)$
the quantized enveloping algebra for the Lie algebra $\lie{sl}_2$.
The category of finite dimensional representations of
$U_q(\lie{sl}_2)$ has a semisimple subquotient $\FF(U_q(\lie{sl}_2))$
called the {\em fusion category} \cite[Section 3.3]{ba:le}.  The
isomorphism classes of simple objects in $\FF(U_q(\lie{sl}_2))$ have
canonical representatives $V_j$ labelled by half-integers $j \in
[0,(r-2)/2]$.  Given $j_1,j_2 \in [0,(r-2)/2] \cap \Z/2$, the tensor
product $V_{j_1} \otimes V_{j_2}$ is isomorphic to the direct sum of
objects $V_{j_3}$ where the sum is over $j_3$ satisfying the {\em
quantum Clebsch-Gordan inequalities}
\begin{equation} \label{qcg} 
\max(j_1 - j_2,j_2-j_1) \leq j_3 \leq  \min(j_1 + j_2,r-2 - j_1 - j_2)
\end{equation} 
and the parity condition $j_1 + j_2 + j_3 \in \Z$.  Geometrically, the
condition \eqref{qcg} means that there exists a triangle in the unit
sphere with edge lengths $j_a(r-2)/2\pi, a = 1,2,3$.  Generalizations
of these inequalities to Lie algebras of higher rank are described in
\cite{ag:ei},\cite{bl:ip},\cite{bel:qhorn},\cite{tw:pa}.

This paper concerns a generalization of this relationship in a
different direction, namely from triangles to tetrahedra.  The {\em
quantum $6j$ symbol} is a function of a $6$-tuple $j_{ab}, 1 \leq a
\leq b \leq 4$, defined as follows: Since the tensor product of simple
modules is multiplicity-free there is a projective basis for
$\Hom_{U_q(\sl_2)}(V_{j_{14}}, V_{j_{12}} \otimes V_{j_{23}} \otimes
V_{j_{34}})$ associated to each way of parenthesizing, parameterized
by half-integers $j_{13}$ resp. $j_{24}$.  The $6j$ symbols $\bsj
j_{12} & j_{23} & j_{13} \\ j_{34} & j_{14} & j_{24} \esj$ are the
coefficients in the change of basis matrix.  $6j$ symbols for $q = 1$
were introduced as a tool in atomic spectroscopy by Racah
\cite{ra:6j}, and then studied mathematically by Wigner \cite{wi:gr}.
$6j$ symbols for $U_q(\lie{sl}_2)$ were introduced by Kirillov and
Reshetikhin \cite{ki:re}, who used them to generalize the Jones knot
invariant. Turaev and Viro used them to define three-manifold
invariants \cite{tv:st}, or what physicists call quantum gravity with
cosmological constant \cite{ba:tqft}.  Roughly speaking these
invariants are obtained by combinatorial integration of a product of
$6j$ symbols: for the Turaev-Viro invariant, a $6j$ symbol is attached
to each tetrahedron in a triangulation, while in the Jones invariant a
$6j$ symbol is attached to each crossing.  Various outstanding
conjectures concern the asymptotics of these invariants: The Witten
conjecture relates the asymptotics of the Turaev-Viro invariant (the
norm-square of the Reshetikhin-Turaev invariant) with Chern-Simons
invariants of flat $SO(4)$ bundles \cite{freed:comp} while the volume
conjecture of Murakami-Murakami relates the asymptotics of the colored
Jones polynomial with the hyperbolic volume of the knot complement
\cite{mu:co}.  A natural approach to these conjectures is to show that
non-Euclidean geometry appears in the asymptotics of the $6j$ symbol;
we will show that this is indeed the case.

The connection of the $6j$ symbols to geometry arises as follows.  By
a theorem of Finkelberg \cite{fi:fu}, $\FF(U_q(\lie{sl}_2))$
is isomorphic to the tensor category of level $r-2$ representations of
the affine Lie algebra $\widehat{\lie{sl}_2}$.  The product for the
latter category uses as its definition the space of genus zero
conformal blocks for Wess-Zumino-Witten (WZW) conformal field theory.
A picture is perhaps the best way of getting across the idea of the
role the quantum $6j$ symbols play in WZW:
\vskip -.4in
\begin{equation} \label{stringy} 
\setlength{\unitlength}{0.00033333in}
\begingroup\makeatletter\ifx\SetFigFont\undefined%
\gdef\SetFigFont#1#2#3#4#5{%
  \reset@font\fontsize{#1}{#2pt}%
  \fontfamily{#3}\fontseries{#4}\fontshape{#5}%
  \selectfont}%
\fi\endgroup%
{\renewcommand{\dashlinestretch}{30}
\begin{picture}(3024,2192)(0,1000)
\put(2112,983){\ellipse{600}{150}}
\put(1512,1883){\ellipse{600}{150}}
\put(312,1883){\ellipse{600}{150}}
\put(2712,1883){\ellipse{600}{150}}
\put(1512,83){\ellipse{600}{150}}
\path(3012,1883)(3011,1882)(3009,1879)
	(3005,1874)(2999,1866)(2991,1855)
	(2981,1842)(2969,1825)(2955,1806)
	(2939,1785)(2922,1762)(2904,1738)
	(2885,1712)(2866,1685)(2846,1658)
	(2826,1630)(2804,1600)(2783,1570)
	(2760,1538)(2737,1505)(2713,1470)
	(2688,1433)(2662,1396)(2637,1358)
	(2608,1314)(2581,1273)(2558,1237)
	(2537,1204)(2519,1174)(2503,1147)
	(2489,1123)(2476,1101)(2465,1081)
	(2454,1063)(2445,1046)(2437,1031)
	(2430,1017)(2424,1006)(2419,997)
	(2416,991)(2414,986)(2413,984)(2412,983)
\path(2412,1883)(2411,1882)(2407,1879)
	(2402,1873)(2394,1865)(2383,1856)
	(2369,1844)(2354,1831)(2337,1817)
	(2319,1804)(2299,1790)(2278,1778)
	(2256,1766)(2232,1756)(2205,1747)
	(2176,1740)(2145,1735)(2112,1733)
	(2079,1735)(2048,1740)(2019,1747)
	(1992,1756)(1968,1766)(1946,1778)
	(1925,1790)(1905,1804)(1887,1817)
	(1870,1831)(1855,1844)(1841,1856)
	(1830,1865)(1822,1873)(1817,1879)
	(1813,1882)(1812,1883)
\path(1212,1883)(1213,1882)(1215,1880)
	(1219,1876)(1225,1870)(1233,1862)
	(1243,1852)(1255,1839)(1269,1824)
	(1285,1808)(1302,1790)(1320,1770)
	(1339,1750)(1358,1728)(1378,1706)
	(1398,1682)(1420,1657)(1441,1631)
	(1464,1603)(1487,1573)(1511,1541)
	(1536,1506)(1562,1470)(1587,1433)
	(1613,1392)(1638,1353)(1660,1316)
	(1680,1282)(1697,1250)(1713,1221)
	(1726,1193)(1739,1166)(1750,1142)
	(1760,1118)(1770,1095)(1778,1074)
	(1786,1055)(1792,1038)(1798,1022)
	(1803,1009)(1806,999)(1809,992)
	(1811,987)(1812,984)(1812,983)
\path(2412,983)(2411,982)(2409,979)
	(2405,974)(2399,966)(2391,955)
	(2381,942)(2369,925)(2355,906)
	(2339,885)(2322,862)(2304,838)
	(2285,812)(2266,785)(2246,758)
	(2226,730)(2204,700)(2183,670)
	(2160,638)(2137,605)(2113,570)
	(2088,533)(2062,496)(2037,458)
	(2008,414)(1981,373)(1958,337)
	(1937,304)(1919,274)(1903,247)
	(1889,223)(1876,201)(1865,181)
	(1854,163)(1845,146)(1837,131)
	(1830,117)(1824,106)(1819,97)
	(1816,91)(1814,86)(1813,84)(1812,83)
\path(12,1883)(13,1882)(14,1880)
	(16,1877)(20,1872)(25,1864)
	(31,1854)(40,1841)(50,1826)
	(62,1808)(75,1788)(90,1766)
	(107,1741)(124,1715)(143,1687)
	(162,1658)(182,1627)(203,1596)
	(225,1563)(248,1529)(272,1493)
	(296,1457)(322,1419)(348,1379)
	(376,1337)(405,1293)(436,1247)
	(469,1198)(503,1147)(538,1094)
	(575,1039)(612,983)(649,927)
	(686,872)(721,819)(755,768)
	(788,719)(819,673)(848,629)
	(876,587)(902,547)(928,509)
	(952,473)(976,437)(999,403)
	(1021,370)(1042,339)(1062,308)
	(1081,279)(1100,251)(1117,225)
	(1134,200)(1149,178)(1162,158)
	(1174,140)(1184,125)(1193,112)
	(1199,102)(1204,94)(1208,89)
	(1210,86)(1211,84)(1212,83)
\path(612,1883)(613,1882)(614,1879)
	(617,1874)(621,1866)(628,1856)
	(636,1841)(646,1823)(659,1802)
	(673,1777)(689,1749)(707,1719)
	(726,1686)(746,1651)(768,1614)
	(790,1576)(813,1537)(837,1497)
	(861,1457)(885,1417)(910,1376)
	(935,1334)(960,1293)(986,1251)
	(1013,1208)(1040,1165)(1068,1122)
	(1096,1079)(1125,1035)(1154,992)
	(1183,949)(1212,908)(1248,857)
	(1281,813)(1309,774)(1332,742)
	(1351,716)(1365,696)(1375,680)
	(1382,668)(1386,659)(1387,653)
	(1388,649)(1387,645)(1386,643)
	(1386,641)(1387,638)(1390,636)
	(1394,632)(1402,628)(1412,623)
	(1426,618)(1444,613)(1464,609)
	(1487,607)(1512,608)(1542,615)
	(1570,627)(1597,644)(1621,664)
	(1644,687)(1665,712)(1685,739)
	(1703,767)(1721,796)(1737,826)
	(1752,855)(1766,884)(1779,910)
	(1789,932)(1798,951)(1804,965)
	(1809,975)(1811,980)(1812,983)
\put(162,2033){\makebox(0,0)[lb]{\smash{{{\SetFigFont{12}{14.4}{\rmdefault}{\mddefault}{\updefault}$j_{12}$}}}}}
\put(1362,2033){\makebox(0,0)[lb]{\smash{{{\SetFigFont{12}{14.4}{\rmdefault}{\mddefault}{\updefault}$j_{23}$}}}}}
\put(2562,2033){\makebox(0,0)[lb]{\smash{{{\SetFigFont{12}{14.4}{\rmdefault}{\mddefault}{\updefault}$j_{34}$}}}}}
\put(1962,1133){\makebox(0,0)[lb]{\smash{{{\SetFigFont{12}{14.4}{\rmdefault}{\mddefault}{\updefault}$j_{24}$}}}}}
\put(1362,233){\makebox(0,0)[lb]{\smash{{{\SetFigFont{12}{14.4}{\rmdefault}{\mddefault}{\updefault}$j_{14}$}}}}}
\end{picture}
}
= \sum_{j_{13}}
 \bsj j_{12} & j_{23} & j_{13} \\ j_{34} & j_{14} &
j_{24} \esj 
\setlength{\unitlength}{0.00033333in}
\begingroup\makeatletter\ifx\SetFigFont\undefined%
\gdef\SetFigFont#1#2#3#4#5{%
  \reset@font\fontsize{#1}{#2pt}%
  \fontfamily{#3}\fontseries{#4}\fontshape{#5}%
  \selectfont}%
\fi\endgroup%
{\renewcommand{\dashlinestretch}{30}
\begin{picture}(3024,2117)(0,1000)
\put(1512,1883){\ellipse{600}{150}}
\put(2712,1883){\ellipse{600}{150}}
\put(312,1883){\ellipse{600}{150}}
\put(1512,83){\ellipse{600}{150}}
\put(912,983){\ellipse{600}{150}}
\path(12,1883)(13,1882)(15,1879)
	(19,1874)(25,1866)(33,1855)
	(43,1842)(55,1825)(69,1806)
	(85,1785)(102,1762)(120,1738)
	(139,1712)(158,1685)(178,1658)
	(198,1630)(220,1600)(241,1570)
	(264,1538)(287,1505)(311,1470)
	(336,1433)(362,1396)(387,1358)
	(416,1314)(443,1273)(466,1237)
	(487,1204)(505,1174)(521,1147)
	(535,1123)(548,1101)(559,1081)
	(570,1063)(579,1046)(587,1031)
	(594,1017)(600,1006)(605,997)
	(608,991)(610,986)(611,984)(612,983)
\path(612,1883)(613,1882)(617,1879)
	(622,1873)(630,1865)(641,1856)
	(655,1844)(670,1831)(687,1817)
	(705,1804)(725,1790)(746,1778)
	(768,1766)(792,1756)(819,1747)
	(848,1740)(879,1735)(912,1733)
	(945,1735)(976,1740)(1005,1747)
	(1032,1756)(1056,1766)(1078,1778)
	(1099,1790)(1119,1804)(1137,1817)
	(1154,1831)(1169,1844)(1183,1856)
	(1194,1865)(1202,1873)(1207,1879)
	(1211,1882)(1212,1883)
\path(1812,1883)(1811,1882)(1809,1880)
	(1805,1876)(1799,1870)(1791,1862)
	(1781,1852)(1769,1839)(1755,1824)
	(1739,1808)(1722,1790)(1704,1770)
	(1685,1750)(1666,1728)(1646,1706)
	(1626,1682)(1604,1657)(1583,1631)
	(1560,1603)(1537,1573)(1513,1541)
	(1488,1506)(1462,1470)(1437,1433)
	(1411,1392)(1386,1353)(1364,1316)
	(1344,1282)(1327,1250)(1311,1221)
	(1298,1193)(1285,1166)(1274,1142)
	(1264,1118)(1254,1095)(1246,1074)
	(1238,1055)(1232,1038)(1226,1022)
	(1221,1009)(1218,999)(1215,992)
	(1213,987)(1212,984)(1212,983)
\path(612,983)(613,982)(615,979)
	(619,974)(625,966)(633,955)
	(643,942)(655,925)(669,906)
	(685,885)(702,862)(720,838)
	(739,812)(758,785)(778,758)
	(798,730)(820,700)(841,670)
	(864,638)(887,605)(911,570)
	(936,533)(962,496)(987,458)
	(1016,414)(1043,373)(1066,337)
	(1087,304)(1105,274)(1121,247)
	(1135,223)(1148,201)(1159,181)
	(1170,163)(1179,146)(1187,131)
	(1194,117)(1200,106)(1205,97)
	(1208,91)(1210,86)(1211,84)(1212,83)
\path(3012,1883)(3011,1882)(3010,1880)
	(3008,1877)(3004,1872)(2999,1864)
	(2993,1854)(2984,1841)(2974,1826)
	(2962,1808)(2949,1788)(2934,1766)
	(2917,1741)(2900,1715)(2881,1687)
	(2862,1658)(2842,1627)(2821,1596)
	(2799,1563)(2776,1529)(2752,1493)
	(2728,1457)(2702,1419)(2676,1379)
	(2648,1337)(2619,1293)(2588,1247)
	(2555,1198)(2521,1147)(2486,1094)
	(2449,1039)(2412,983)(2375,927)
	(2338,872)(2303,819)(2269,768)
	(2236,719)(2205,673)(2176,629)
	(2148,587)(2122,547)(2096,509)
	(2072,473)(2048,437)(2025,403)
	(2003,370)(1982,339)(1962,308)
	(1943,279)(1924,251)(1907,225)
	(1890,200)(1875,178)(1862,158)
	(1850,140)(1840,125)(1831,112)
	(1825,102)(1820,94)(1816,89)
	(1814,86)(1813,84)(1812,83)
\path(2412,1883)(2411,1882)(2410,1880)
	(2407,1875)(2403,1868)(2396,1859)
	(2388,1846)(2378,1830)(2365,1811)
	(2351,1789)(2335,1764)(2317,1737)
	(2298,1707)(2278,1676)(2256,1643)
	(2234,1609)(2211,1574)(2187,1538)
	(2163,1501)(2139,1464)(2114,1427)
	(2089,1389)(2064,1351)(2038,1312)
	(2011,1272)(1984,1232)(1956,1191)
	(1928,1150)(1899,1108)(1870,1065)
	(1841,1024)(1812,983)(1773,928)
	(1738,878)(1708,836)(1684,800)
	(1666,771)(1653,748)(1644,730)
	(1639,716)(1637,705)(1636,696)
	(1637,689)(1638,683)(1638,677)
	(1636,670)(1633,663)(1627,655)
	(1617,646)(1603,636)(1586,627)
	(1564,618)(1539,611)(1512,608)
	(1482,612)(1454,621)(1427,636)
	(1403,655)(1380,677)(1359,702)
	(1339,729)(1321,758)(1303,788)
	(1287,819)(1272,849)(1258,879)
	(1245,906)(1235,930)(1226,950)
	(1220,965)(1215,975)(1213,980)(1212,983)
\put(162,1958){\makebox(0,0)[lb]{\smash{{{\SetFigFont{12}{14.4}{\rmdefault}{\mddefault}{\updefault}$j_{12}$}}}}}
\put(1437,1958){\makebox(0,0)[lb]{\smash{{{\SetFigFont{12}{14.4}{\rmdefault}{\mddefault}{\updefault}$j_{23}$}}}}}
\put(2562,1958){\makebox(0,0)[lb]{\smash{{{\SetFigFont{12}{14.4}{\rmdefault}{\mddefault}{\updefault}$j_{34}$}}}}}
\put(687,1058){\makebox(0,0)[lb]{\smash{{{\SetFigFont{12}{14.4}{\rmdefault}{\mddefault}{\updefault}$j_{13}$}}}}}
\put(1362,158){\makebox(0,0)[lb]{\smash{{{\SetFigFont{12}{14.4}{\rmdefault}{\mddefault}{\updefault}$j_{14}$}}}}}
\end{picture}
}.
\end{equation}
\vskip .2in
\noindent The vector space of conformal blocks is isomorphic to the
space of holomorphic sections of the determinant line bundle on the
moduli space of flat $SU(2)$-bundles, see Pauly \cite{pa:co}.  The
$6j$ symbol is equal to the Hermitian pairing between two holomorphic
sections.  According to geometric quantization, the sections should be
interpreted as quantum\footnote{In this subject the term {\em quantum}
is used rather confusingly in two senses: $ q \neq 1 $ and $\hbar \neq
0$.  From our point of view, $-\ln(q)^2$ describes the curvature or,
in physics language, the cosmological constant.  The true quantum
parameter is $\hbar$, which appears in the guise of $1/k$ where
quantum states are sections of the $k$-th tensor power of the
determinant bundle.}-mechanical states, and one should try to show
that in the semiclassical limit the states concentrate to Lagrangian
submanifolds.  Assuming this holds, the leading term in the
asymptotics of the pairing is a sum over intersection points of the
Lagrangians, involving Poisson brackets of the functions defining the
Lagrangians and the pairing between the sections in the fiber of the
determinant line bundle, see Borthwick-Paul-Uribe \cite{bo:le}.  The
particular case of interest is the moduli space of flat $SU(2)$
bundles on the four-holed two-sphere, with holonomies around the four
boundary components fixed by the labels $j_{12},j_{23},j_{34},j_{14}$.
The Lagrangians consists of flat bundles for which the holonomy around
the intermediate circle correspond to $j_{13},j_{24}$ respectively, as
in \eqref{stringy}.  Using the diffeomorphism $SU(2) \to S^3$, the
intersection points correspond to spherical tetrahedra with edge
lengths $l_{ab} := \pi j_{ab}/(r-2)$.
\begin{figure}[h]
\setlength{\unitlength}{0.00033333in}
%
{\renewcommand{\dashlinestretch}{30}
\begin{picture}(7524,3639)(-2000,-10)
\path(1512,3612)(1513,3612)(1515,3611)
	(1518,3610)(1523,3609)(1531,3607)
	(1541,3604)(1554,3600)(1569,3596)
	(1588,3590)(1610,3584)(1635,3577)
	(1662,3569)(1692,3560)(1724,3550)
	(1759,3540)(1795,3528)(1834,3516)
	(1874,3504)(1915,3490)(1958,3476)
	(2002,3461)(2046,3446)(2092,3430)
	(2139,3413)(2187,3395)(2236,3376)
	(2286,3357)(2338,3336)(2391,3314)
	(2445,3291)(2502,3266)(2560,3240)
	(2619,3213)(2681,3183)(2745,3152)
	(2810,3119)(2877,3085)(2944,3049)
	(3012,3012)(3082,2972)(3151,2932)
	(3217,2893)(3281,2853)(3341,2815)
	(3399,2778)(3453,2741)(3505,2706)
	(3554,2671)(3601,2638)(3645,2605)
	(3688,2573)(3728,2541)(3767,2511)
	(3804,2480)(3840,2451)(3875,2421)
	(3908,2393)(3940,2365)(3971,2338)
	(4000,2312)(4028,2287)(4054,2263)
	(4078,2240)(4100,2219)(4121,2200)
	(4139,2183)(4155,2167)(4169,2154)
	(4181,2142)(4191,2133)(4198,2126)
	(4204,2120)(4208,2116)(4210,2114)
	(4211,2113)(4212,2112)
\path(4212,2112)(4211,2112)(4210,2111)
	(4206,2109)(4201,2106)(4194,2102)
	(4184,2097)(4172,2090)(4157,2081)
	(4139,2071)(4118,2060)(4095,2046)
	(4070,2032)(4042,2015)(4012,1998)
	(3981,1979)(3947,1959)(3913,1937)
	(3877,1915)(3840,1892)(3802,1867)
	(3763,1842)(3724,1815)(3684,1788)
	(3643,1759)(3601,1728)(3558,1697)
	(3514,1663)(3469,1628)(3423,1591)
	(3376,1552)(3327,1510)(3277,1466)
	(3225,1420)(3173,1371)(3120,1320)
	(3066,1267)(3012,1212)(2959,1156)
	(2908,1101)(2860,1046)(2813,992)
	(2769,940)(2728,889)(2689,840)
	(2652,792)(2618,745)(2586,700)
	(2555,656)(2526,614)(2498,572)
	(2472,531)(2447,491)(2423,452)
	(2401,414)(2379,376)(2358,340)
	(2338,305)(2320,271)(2302,238)
	(2286,207)(2270,178)(2256,151)
	(2243,126)(2232,103)(2222,83)
	(2213,66)(2206,51)(2200,39)
	(2196,30)(2192,23)(2190,18)
	(2188,14)(2187,13)(2187,12)
\path(1512,3612)(1511,3612)(1509,3610)
	(1506,3609)(1501,3605)(1493,3601)
	(1483,3595)(1471,3587)(1456,3578)
	(1438,3567)(1418,3554)(1396,3540)
	(1372,3523)(1346,3506)(1318,3486)
	(1289,3466)(1259,3443)(1228,3419)
	(1196,3394)(1163,3367)(1130,3339)
	(1095,3308)(1060,3276)(1025,3241)
	(988,3204)(951,3164)(912,3120)
	(872,3074)(831,3023)(789,2969)
	(745,2910)(701,2848)(656,2781)
	(612,2712)(575,2651)(539,2590)
	(505,2529)(472,2469)(441,2409)
	(411,2350)(384,2293)(358,2237)
	(333,2183)(310,2129)(288,2077)
	(267,2025)(248,1974)(229,1925)
	(211,1876)(194,1828)(178,1780)
	(163,1734)(148,1688)(134,1643)
	(120,1600)(108,1557)(96,1517)
	(84,1478)(74,1441)(64,1406)
	(55,1374)(47,1344)(40,1318)
	(34,1295)(28,1274)(24,1257)
	(20,1243)(17,1232)(15,1224)
	(14,1218)(13,1215)(12,1213)(12,1212)
\path(12,1212)(13,1212)(15,1212)
	(19,1212)(25,1212)(34,1212)
	(46,1212)(61,1211)(79,1211)
	(101,1210)(125,1209)(153,1208)
	(183,1207)(216,1205)(251,1203)
	(287,1200)(326,1197)(366,1193)
	(407,1189)(450,1184)(493,1178)
	(537,1171)(582,1163)(629,1155)
	(676,1145)(725,1133)(774,1120)
	(825,1106)(878,1089)(932,1071)
	(988,1050)(1046,1028)(1105,1002)
	(1165,974)(1226,944)(1287,912)
	(1350,876)(1411,838)(1469,800)
	(1524,762)(1575,724)(1623,686)
	(1668,648)(1710,611)(1750,574)
	(1787,538)(1822,502)(1855,466)
	(1886,431)(1916,396)(1945,361)
	(1972,327)(1997,293)(2021,260)
	(2044,229)(2066,198)(2086,169)
	(2104,142)(2120,118)(2135,95)
	(2147,75)(2158,59)(2167,44)
	(2174,33)(2179,25)(2183,19)
	(2185,15)(2186,13)(2187,12)
\dashline{60.000}(12,1212)(13,1212)(14,1214)
	(16,1216)(20,1219)(25,1224)
	(33,1230)(42,1238)(54,1248)
	(68,1259)(85,1273)(104,1289)
	(125,1306)(149,1325)(175,1345)
	(203,1367)(233,1391)(266,1415)
	(300,1441)(336,1467)(374,1494)
	(413,1521)(454,1549)(496,1577)
	(540,1606)(585,1634)(632,1662)
	(681,1691)(731,1719)(783,1747)
	(838,1775)(895,1803)(954,1831)
	(1016,1859)(1081,1886)(1149,1913)
	(1220,1940)(1294,1967)(1373,1993)
	(1454,2019)(1540,2044)(1628,2068)
	(1719,2091)(1812,2112)(1898,2130)
	(1985,2146)(2070,2161)(2155,2174)
	(2238,2185)(2319,2195)(2399,2204)
	(2476,2211)(2552,2216)(2625,2221)
	(2697,2224)(2767,2227)(2836,2228)
	(2903,2228)(2968,2228)(3033,2227)
	(3096,2225)(3158,2222)(3219,2219)
	(3278,2215)(3337,2211)(3395,2207)
	(3451,2202)(3507,2196)(3560,2191)
	(3613,2185)(3664,2179)(3713,2173)
	(3760,2167)(3805,2162)(3847,2156)
	(3887,2150)(3924,2145)(3958,2140)
	(3990,2135)(4018,2131)(4042,2127)
	(4064,2124)(4082,2121)(4098,2119)
	(4110,2117)(4120,2115)(4127,2114)
	(4132,2113)(4135,2112)(4136,2112)(4137,2112)
\dashline{60.000}(1512,3612)(1512,3611)(1513,3610)
	(1515,3606)(1517,3601)(1520,3594)
	(1525,3584)(1531,3571)(1538,3556)
	(1546,3537)(1556,3515)(1567,3491)
	(1580,3463)(1593,3433)(1608,3400)
	(1623,3364)(1640,3326)(1657,3286)
	(1674,3244)(1693,3200)(1711,3155)
	(1730,3108)(1749,3060)(1769,3011)
	(1788,2960)(1808,2908)(1827,2854)
	(1847,2799)(1866,2743)(1886,2685)
	(1905,2625)(1925,2562)(1944,2498)
	(1964,2432)(1983,2363)(2003,2291)
	(2022,2217)(2041,2140)(2060,2060)
	(2078,1979)(2095,1896)(2112,1812)
	(2128,1724)(2143,1638)(2156,1553)
	(2168,1471)(2179,1393)(2188,1317)
	(2195,1244)(2202,1175)(2207,1108)
	(2212,1044)(2215,983)(2218,923)
	(2220,866)(2221,810)(2222,757)
	(2222,705)(2221,654)(2220,605)
	(2219,557)(2217,511)(2216,467)
	(2213,424)(2211,383)(2209,344)
	(2206,308)(2204,274)(2201,242)
	(2199,214)(2197,188)(2195,166)
	(2193,146)(2191,130)(2190,117)
	(2189,106)(2188,98)(2188,93)
	(2187,90)(2187,88)(2187,87)
\put(237,2712){\makebox(0,0)[lb]{\smash{{{\SetFigFont{12}{14.4}{\rmdefault}{\mddefault}{\updefault}$l_{12}$}}}}}
\put(2787,3237){\makebox(0,0)[lb]{\smash{{{\SetFigFont{12}{14.4}{\rmdefault}{\mddefault}{\updefault}$l_{23}$}}}}}
\put(3162,912){\makebox(0,0)[lb]{\smash{{{\SetFigFont{12}{14.4}{\rmdefault}{\mddefault}{\updefault}$l_{34}$}}}}}
\put(837,612){\makebox(0,0)[lb]{\smash{{{\SetFigFont{12}{14.4}{\rmdefault}{\mddefault}{\updefault}$l_{14}$}}}}}
\put(1962,2712){\makebox(0,0)[lb]{\smash{{{\SetFigFont{12}{14.4}{\rmdefault}{\mddefault}{\updefault}$l_{13}$}}}}}
\put(1062,1587){\makebox(0,0)[lb]{\smash{{{\SetFigFont{12}{14.4}{\rmdefault}{\mddefault}{\updefault}$l_{24}$}}}}}
\end{picture}
}
\caption{The spherical tetrahedon $\tau$}
\end{figure}
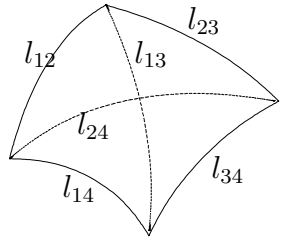
\noindent 
The Poisson bracket turns out to be the determinant of the Gram matrix
of the tetrahedron, and the phase shift is the area of a holomorphic
disk in the moduli space which can be computed using Schl\"afli's
formula.  These ideas lead to the
conjecture that if $\tau$ is non-degenerate then
\begin{equation}  \label{conj}
\bsj kj_{12} & kj_{23} & kj_{13} \\ kj_{34} & kj_{14} & kj_{24}
  \esj_{q = \exp(\pi i /r(k))} \sim \frac{2\pi \cos( \phi(k) +
  \pi/4)}{r(k)^{3/2} \det(\cos(l_{ab}))^{1/4}}\end{equation}
as $k \to \infty$ where $r(k) = k(r-2) + 2$, 
$$\phi(k) = \frac{r(k)}{2\pi} \left(\sum_{a < b} \theta_{ab}(k) l_{ab}(k) -
2\Vol(\tau(k)) \right)$$
$\tau(k)$ is the spherical tetrahedron with edge lengths
$$ l_{ab}(k) = 2\pi \frac{kj_{ab} + \hh}{r(k)} ,$$ 
and $\theta_{ab}(k)$ are its exterior dihedral angles.  The formula
\eqref{conj} generalizes one for the case $q = 1$ given by the
physicists Wigner \cite{wi:gr}, Ponzano and Regge \cite{po:6j}, and
proved by J. Roberts \cite{ro:qu}, see also \cite{ch:as},\cite{wi:as}.
The numerator of the formula was conjectured by Mizoguchi and Tada
\cite{mi:th}.  See also Roberts' discussion \cite{ro:as}.  

Our proof of \eqref{conj} is not based on the above ideas.  Instead,
we show (following Schulten-Gordon \cite{sch:semi} and \cite{mi:th})
that both sides satisfy a second order difference equation as one
label is varied.  It follows that each side is a linear combination of
the two linearly independent solutions to the equation.  Taking the
Euclidean limit and applying Roberts' theorem \cite{ro:qu} shows
equality of the coefficients.  This method, while less geometric, has
the advantage that it works for non-integral $r$ and gives a formula
for degenerate tetrahedra with non-degenerate faces.  In fact there
are seven different cases depending on the type of degeneracy of the
tetrahedron $\tau$; these are pictured in Figure 2.  The second
column, titled Maslov diagram, shows the geometry of the relevant
Lagrangian submanifolds (or, for the more degenerate cases, isotropic
submanifolds).
\begin{figure}[h]
\setlength{\unitlength}{0.00023333in}
\begingroup\makeatletter\ifx\SetFigFont\undefined%
\gdef\SetFigFont#1#2#3#4#5{%
  \reset@font\fontsize{#1}{#2pt}%
  \fontfamily{#3}\fontseries{#4}\fontshape{#5}%
  \selectfont}%
\fi\endgroup%
\renewcommand{\dashlinestretch}{30}
\begin{picture}(7035,17049)(0,-10)
\put(4605,14190){\ellipse{618}{618}}
\put(4605,14190){\ellipse{1210}{1210}}
\put(4605,14190){\ellipse{1806}{1806}}
\put(4605,14190){\ellipse{2418}{2418}}
\put(5805,14190){\ellipse{600}{600}}
\put(5805,14190){\ellipse{1200}{1200}}
\put(5805,14190){\ellipse{1800}{1800}}
\put(5805,14190){\ellipse{2400}{2400}}
\put(4605,7290){\ellipse{618}{618}}
\put(4605,7290){\ellipse{1210}{1210}}
\put(4605,7290){\ellipse{1806}{1806}}
\put(4605,7290){\ellipse{2418}{2418}}
\put(5805,7290){\ellipse{600}{600}}
\put(5805,7290){\ellipse{1200}{1200}}
\put(5805,7290){\ellipse{1800}{1800}}
\put(5805,7290){\ellipse{2400}{2400}}
\put(4605,10590){\ellipse{618}{618}}
\put(4605,10590){\ellipse{1210}{1210}}
\put(4605,10590){\ellipse{1806}{1806}}
\put(4605,10590){\ellipse{2418}{2418}}
\put(5805,10590){\ellipse{600}{600}}
\put(5805,10590){\ellipse{1200}{1200}}
\put(5805,10590){\ellipse{1800}{1800}}
\put(5805,10590){\ellipse{2400}{2400}} \texture{44555555 55aaaaaa
aa555555 55aaaaaa aa555555 55aaaaaa aa555555 55aaaaaa aa555555
55aaaaaa aa555555 55aaaaaa aa555555 55aaaaaa aa555555 55aaaaaa
aa555555 55aaaaaa aa555555 55aaaaaa aa555555 55aaaaaa aa555555
55aaaaaa aa555555 55aaaaaa aa555555 55aaaaaa aa555555 55aaaaaa
aa555555 55aaaaaa } \put(705,15390){\shade\ellipse{150}{150}}
\put(705,15390){\ellipse{150}{150}}
\put(2505,14790){\shade\ellipse{150}{150}}
\put(2505,14790){\ellipse{150}{150}}
\put(1305,12990){\shade\ellipse{150}{150}}
\put(1305,12990){\ellipse{150}{150}}
\put(105,14190){\shade\ellipse{150}{150}}
\put(105,14190){\ellipse{150}{150}}
\put(705,11790){\shade\ellipse{150}{150}}
\put(705,11790){\ellipse{150}{150}}
\put(2505,11190){\shade\ellipse{150}{150}}
\put(2505,11190){\ellipse{150}{150}}
\put(1305,9390){\shade\ellipse{150}{150}}
\put(1305,9390){\ellipse{150}{150}}
\put(105,10590){\shade\ellipse{150}{150}}
\put(105,10590){\ellipse{150}{150}}
\put(105,8190){\shade\ellipse{150}{150}}
\put(105,8190){\ellipse{150}{150}}
\put(2505,7590){\shade\ellipse{150}{150}}
\put(2505,7590){\ellipse{150}{150}}
\put(105,6390){\shade\ellipse{150}{150}}
\put(105,6390){\ellipse{150}{150}}
\put(105,7590){\shade\ellipse{150}{150}}
\put(105,7590){\ellipse{150}{150}} \thicklines
\put(105,5190){\shade\ellipse{150}{150}}
\put(105,5190){\ellipse{150}{150}}
\put(2505,4590){\shade\ellipse{150}{150}}
\put(2505,4590){\ellipse{150}{150}}
\put(180,3390){\shade\ellipse{150}{150}}
\put(180,3390){\ellipse{150}{150}}
\put(4605,14190){\ellipse{1800}{1800}}
\put(5805,14190){\ellipse{1800}{1800}}
\put(4605,10590){\ellipse{1200}{1200}}
\put(5805,10590){\ellipse{1200}{1200}} \texture{55888888 88555555
5522a222 a2555555 55888888 88555555 552a2a2a 2a555555 55888888
88555555 55a222a2 22555555 55888888 88555555 552a2a2a 2a555555
55888888 88555555 5522a222 a2555555 55888888 88555555 552a2a2a
2a555555 55888888 88555555 55a222a2 22555555 55888888 88555555
552a2a2a 2a555555 } \put(4605,7290){\shade\ellipse{150}{150}}
\put(4605,7290){\ellipse{150}{150}}
\put(5805,7290){\ellipse{2404}{2404}} \texture{44555555 55aaaaaa
aa555555 55aaaaaa aa555555 55aaaaaa aa555555 55aaaaaa aa555555
55aaaaaa aa555555 55aaaaaa aa555555 55aaaaaa aa555555 55aaaaaa
aa555555 55aaaaaa aa555555 55aaaaaa aa555555 55aaaaaa aa555555
55aaaaaa aa555555 55aaaaaa aa555555 55aaaaaa aa555555 55aaaaaa
aa555555 55aaaaaa } \put(2535,2235){\shade\ellipse{150}{150}}
\put(2535,2235){\ellipse{150}{150}}
\put(705,2220){\shade\ellipse{150}{150}}
\put(705,2220){\ellipse{150}{150}}
\put(90,2220){\shade\ellipse{150}{150}}
\put(90,2220){\ellipse{150}{150}}
\put(1605,2220){\shade\ellipse{150}{150}}
\put(1605,2220){\ellipse{150}{150}}
\put(105,990){\shade\ellipse{150}{150}}
\put(105,990){\ellipse{150}{150}}
\put(2505,990){\shade\ellipse{150}{150}}
\put(2505,990){\ellipse{150}{150}} \texture{55888888 88555555 5522a222
a2555555 55888888 88555555 552a2a2a 2a555555 55888888 88555555
55a222a2 22555555 55888888 88555555 552a2a2a 2a555555 55888888
88555555 5522a222 a2555555 55888888 88555555 552a2a2a 2a555555
55888888 88555555 55a222a2 22555555 55888888 88555555 552a2a2a
2a555555 } \put(4905,990){\shade\ellipse{150}{150}}
\put(4905,990){\ellipse{150}{150}} \texture{44555555 55aaaaaa aa555555
55aaaaaa aa555555 55aaaaaa aa555555 55aaaaaa aa555555 55aaaaaa
aa555555 55aaaaaa aa555555 55aaaaaa aa555555 55aaaaaa aa555555
55aaaaaa aa555555 55aaaaaa aa555555 55aaaaaa aa555555 55aaaaaa
aa555555 55aaaaaa aa555555 55aaaaaa aa555555 55aaaaaa aa555555
55aaaaaa } \put(105,90){\shade\ellipse{150}{150}}
\put(105,90){\ellipse{150}{150}} \texture{55888888 88555555 5522a222
a2555555 55888888 88555555 552a2a2a 2a555555 55888888 88555555
55a222a2 22555555 55888888 88555555 552a2a2a 2a555555 55888888
88555555 5522a222 a2555555 55888888 88555555 552a2a2a 2a555555
55888888 88555555 55a222a2 22555555 55888888 88555555 552a2a2a
2a555555 } \put(4905,90){\shade\ellipse{150}{150}}
\put(4905,90){\ellipse{150}{150}}
\put(4905,2205){\shade\ellipse{150}{150}}
\put(4905,2205){\ellipse{150}{150}}\thicklines
\put(6117,4292){\ellipse{1800}{1800}}  \thinlines
\put(6105,4290){\ellipse{1200}{1200}}
\put(6089,4265){\ellipse{600}{600}} \thicklines
\put(4610,4290){\ellipse{1210}{1210}} \thinlines
\put(4608,4290){\ellipse{1806}{1806}}
\put(4605,4290){\ellipse{2404}{2404}} \thicklines
\path(705,15390)(1305,12990) \path(105,14190)(705,14340)
\path(1155,14490)(2505,14790) \path(105,10590)(705,11790)(2505,11190)
(1305,9390)(105,10590)(2505,11190) \path(705,11790)(1305,9390)
\path(105,7590)(2505,7590)
\path(105,8190)(105,6390)(2505,7590)(105,8190) \texture{44555555
55aaaaaa aa555555 55aaaaaa aa555555 55aaaaaa aa555555 55aaaaaa
aa555555 55aaaaaa aa555555 55aaaaaa aa555555 55aaaaaa aa555555
55aaaaaa aa555555 55aaaaaa aa555555 55aaaaaa aa555555 55aaaaaa
aa555555 55aaaaaa aa555555 55aaaaaa aa555555 55aaaaaa aa555555
55aaaaaa aa555555 55aaaaaa } \path(180,5190)(180,3390)
\path(180,5190)(180,3390) \path(109,5133)(2509,4533)
\path(109,5133)(2509,4533) \path(105,3390)(2505,4590)
\path(105,3390)(2505,4590) \path(120,5190)(120,3390)
\path(120,5190)(120,3390) \path(109,5189)(2509,4589)
\path(109,5189)(2509,4589) \path(105,14190)(705,15390)(2505,14790)
(1305,12990)(105,14190) \path(105,2190)(2580,2190)
\path(135,2250)(2610,2250) \path(105,990)(2505,990)
\path(105,990)(2505,990) \path(105,1050)(2505,1050)
\path(105,1050)(2505,1050)
\put(105,15990){\makebox(0,0)[lb]{\smash{{{Tetrahedron}}}}}
\put(3820,16005){\makebox(0,0)[lb]{\smash{{{Maslov Diagram}}}}}
\put(8505,14490){\makebox(0,0)[lb]{\smash{{{${-3/2}$}}}}}
\put(8505,7290){\makebox(0,0)[lb]{\smash{{{${-5/4}$?}}}}}
\put(8505,10590){\makebox(0,0)[lb]{\smash{{{${-4/3}$}}}}}
\put(8505,4290){\makebox(0,0)[lb]{\smash{{{${-1}$}}}}}
\put(8505,990){\makebox(0,0)[lb]{\smash{{{${-1/2}$}}}}}
\put(8505,90){\makebox(0,0)[lb]{\smash{{{$0$}}}}}
\put(8505,2190){\makebox(0,0)[lb]{\smash{{{${-1}$}}}}}
\put(8505,15990){\makebox(0,0)[lb]{\smash{{{Power of $k$}}}}}
\put(-1505,15990){\makebox(0,0)[lb]{\smash{{{Case}}}}}
\put(-1505,14490){\makebox(0,0)[lb]{\smash{{{(a)}}}}}
\put(-1505,10590){\makebox(0,0)[lb]{\smash{{{(b)}}}}}
\put(-1505,7290){\makebox(0,0)[lb]{\smash{{{(c)}}}}}
\put(-1505,4290){\makebox(0,0)[lb]{\smash{{{(d)}}}}}
\put(-1505,2190){\makebox(0,0)[lb]{\smash{{{(e)}}}}}
\put(-1505,990){\makebox(0,0)[lb]{\smash{{{(f)}}}}}
\put(-1505,90){\makebox(0,0)[lb]{\smash{{{(g)}}}}}

\end{picture}

\caption{Seven types of tetrahedra and the corresponding 
Lagrangians}
\end{figure}
The possible degenerations in the table are
$$ (a) \to (b) \to (c) \to ( (d) \text{\ or \ } (e) ) \to (f) \to (g) .$$
We refer to case (a) as the non-degenerate case, and case (b) as the
tangent case.  There is an overlap between cases (d) and (e), when the
vertices are colinear but $l_{23} = 0$; perhaps cases (d) and (e)
should not really be considered separate.  We managed to cover them
almost completely. The exception is the case (c) that exactly one face
is degenerate: although one can compute the asymptotics, we did not
find a simple formula or manage to construct non-trivial cases where
the configuration occurs.  We also failed to cover the eighth
``classically forbidden'' case that the tetrahedron does not exist,
where one expects (and numerically experiments show) exponential decay
of the $6j$ symbols.  The main theorem (case (c) is intentionally
missing) follows:
\begin{theorem} \label{main}  Let $r >2, j_{12},\ldots,j_{34}
\in [0,(r-2)/2] \cap \Z/2$ and
$$s(k) := \bsj kj_{12} & kj_{23} &
 kj_{13} \\ kj_{34} & kj_{14} & kj_{24} \esj_{q = \exp(\pi i /r(k))} .$$
\begin{enumerate}
\item If $\tau$ exists and is non-degenerate, then
%
$$ s(k) \sim \frac{2\pi \cos( \phi(k) + \pi/4)}{r(k)^{3/2}
  \det(\cos(l_{ab}))^{1/4}} .$$
%
\item If $\tau$ exists, is degenerate and all faces are non-degenerate
then
$$
s(k) \sim r(k)^{-\frac{4}{3}} 2^{\frac{2}{3}} 3^{-\frac{1}{3}} 
\pi^{4/3} \Gamma\left(\frac{2}{3}\right)^{-1} 
\frac{\cos(k\sum \theta_{ab} j_{ab}) }{(A_1 A_2
A_3 A_4)^{1/6}} $$
where $ A_a = \det(\cos(l_{bc})_{b,c \neq a})^{1/2} $ and $\cos(k \sum
\theta_{ab} j_{ab}) = \pm 1$ is a sign depending on the sum of the
$j_{ab}$'s for which $\theta_{ab} = \pi$.
%
%
%
%
\item[(d)]   If $\tau$ exists and has exactly one edge length, say
$l_{ab}$ vanishing, then
$$ 
s(k) \sim 
(-1)^{k(j_{bc} + j_{cd} + j_{bd})}
\pi r(k)^{-1} ( \sin(l_{ac})\sin(l_{bd}))^{-1/2} .$$
Note that $\tau$ need not have any non-degenerate faces.
\item[(e)] If $\tau$ exists, all faces are degenerate, the
vertices lie on a geodesic of length at most $\pi$ in order $a,b,c,d$
and the edge lengths $l_{ac},l_{bd}$ are non-vanishing, then
$$ 
s(k) \sim (-1)^{2kj_{ad}}
 \pi r(k)^{-1} (\sin(l_{ac})
\sin(l_{bd}))^{-1/2}.$$
Note the edge length $l_{bc}$ may vanish.
\item[(f)] If $\tau$ exists, all faces are degenerate, either 
\begin{enumerate} 
\item all non-zero edge lengths are equal and less than
  $\pi$, or 
\item one edge length is zero and the opposite edge length is $\pi$
\end{enumerate}
and $l_{ab} \in (0,\pi)$ then 
$$
%
s(k) \sim  
(-1)^{2kj_{ab}} \pi^{1/2} r(k)^{-1/2}  \sin(l_{ab})^{-1/2} .$$
\item[(g)] If $\tau$ exists, but all edge lengths are $0$ or $\pi$
then $ s(k) = 1 .$
\end{enumerate}
\end{theorem}

Some of these formulas hold for slightly more general sequences of
labels, see Theorems \ref{gen1},\ref{gen2} below.  Most of the results
were checked numerically.  In Figure 3 we show a numerical comparison
generated by Maple between the $6j$ symbols 
and the
asymptotic formula Theorem \ref{main} (a) 
for the
sequence of $6j$ symbols
$ \bsj 40 & 48 & 50 \\
        52 & 54 & j \esj $
as $j$ varies from $0$ to $108$, with $r=200$.\footnote{For the
purposes of Maple, dominant weights were labelled by non-negative {\em
integers}.}  The phase function $\phi$ was computed by numerical
integration, using Schl\"afli's formula \ref{nonEuclid}
\eqref{schlafli}.  The amplitudes $ \pm 2\pi r(k)^{-\thh}
\det(\cos(l_{ab}))^{-\qq}$ are also shown
, as well as the functions from case (b) of the main theorem governing the
degenerate limits
.

\begin{figure}[ht]
\[\psfig{figure=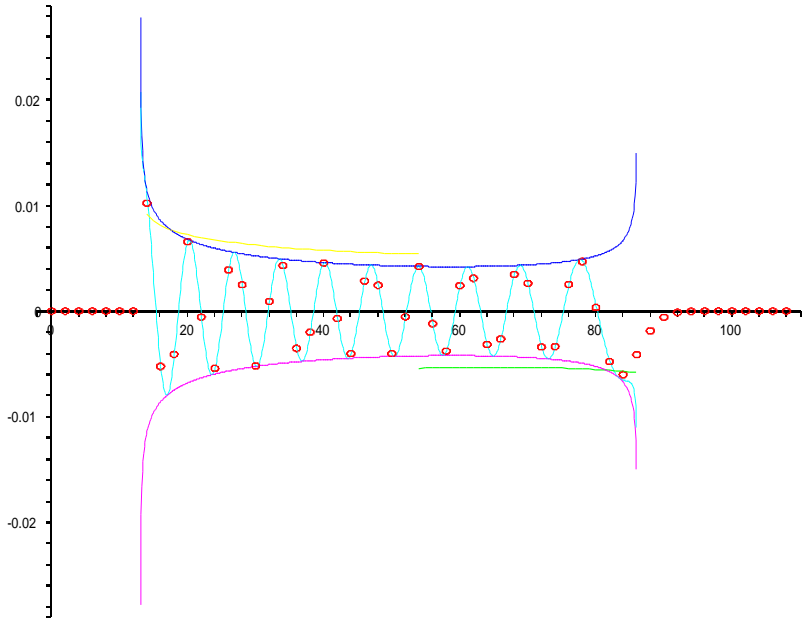,width=5in,height=4in}\]
\caption{Quantum $6j$ symbols versus the asymptotic formulas}
\end{figure}

An outline of the paper is as follows.  Section 2 contains 
background.  The conjecture on the geometry of
conformal blocks in the semiclassical limit is explained in Section 3.
Sections 4 to 7 contain a proof of the main theorem.  Section 8
contains a short discussion of the classically forbidden case. The
results for the hyperbolic case, that is, the case $q$ real and
positive, are stated without proof in Section 9.  Various
questions we were not able to resolve are listed
in Section 10.

{\em Acknowledgements} The project was suggested by J. Roberts at the
end of \cite{ro:qu}, and we thank him for his encouragement.
Discussions with L. Rozansky, F. Luo, I. Korepanov, M. Leingang,
M. McDuffee and I. Rivin were also helpful.  The reference
\cite{sch:semi} and the idea of using recursion to prove the formulas
were pointed out to us by N. Reshetikhin.  Related work on the
asymptotics of Turaev-Viro has been carried out by Frohman and
Kania-Bartoszynska \cite{fro:qu} and J. Murakami, unpublished.
 
\section{Background}

\subsection{The quantum group $U_q(\lie{sl}_2)$}

The quantum group $U_q(\lie{sl}_2)$ was introduced by Kulish and
Reshetikhin in \cite{kr:uq}.  Let $q$ be a complex number not equal to
$0,1$ or $-1$.  For any integer $n$, the {\em quantum integer} $[n]$
is defined by
$$ [n] = \frac{q^n - q^{-n}}{q - q^{-1} } .$$
Note that $[n] \to n$ as $q \to 1$, and also 
$$ [n] =
\frac{ \sin(n\pi/r)}{ \sin(\pi/r)}
           \ \ \ \ \ \      \text{if} \ q = \exp(\pi i/r) .$$
Let $\lie{sl}_2$ denote the Lie algebra of $2\times2$ traceless
matrices, with basis $E,F,H$ and relations 
$ [H,E] = E, \ \ [H,F] = -F, \ \ [E,F] = 2H .$
The irreducible representations of $\sl_2$ are $V_j, \ j =
0,\hh,1,\ldots,$ where $V_j$ is the space of degree $2j$ homogeneous
polynomials in two variables $z,w$.  A basis for $V_j$ is given by
$v_{j,j}, v_{j,j-1}, \ldots, v_{j,-j}$ where the {\em weight vector}
$v_{j,m}$ is defined by 
$$ v_{j,m} = z^a w^b, \ \ a - b = 2m, \ \ a + b = 2 j. $$
The action of the $\sl_2$ in this basis is given by 
$$ Ev_{j,m} = (j-m)v_{j,m+1}, \ \ Fv_{j,m} = (j+m)v_{j,m-1},
\ \ Hv_{j,m} = mv_{j,m} .$$
The quantum group $U_q(\lie{sl}_2)$ is the algebra with generators $E,F,K,K^{-1}$
and relations
$$ [E,F] = \frac{K^2 - K^{-2}}{q - q^{-1}}, \ \ KE = qEK, \ \ \ KF =
q^{-1} FK .$$
The action of $U_q(\lie{sl}_2)$ on $V_j$ is defined by 
$$ Ev_{j,m} = [j-m] v_{j,m+1},
\ \ Fv_{j,m} = [j+m] v_{j,m-1}, \ \ Kv_{j,m} = q^m v_{j,m} .$$
The coproduct on $U_q(\lie{sl}_2)$ defined by 
$$ E \mapsto E \otimes K  +  K^{-1} \otimes E, 
\ \ F \mapsto F \otimes K + K^{-1} \otimes F,
\ \ K^{\pm 1} \mapsto K^{\pm 1} \otimes K^{\pm 1} $$
makes the category $\RR(U_q(\lie{sl}_2))$ of finite dimensional
$U_q(\lie{sl}_2)$ modules a tensor category.  $\RR(U_q(\lie{sl}_2))$
has a semisimple subquotient $\FF(U_q(\lie{sl}_2))$ called the fusion
category for $U_q(\lie{sl}_2)$ \cite[Section 3.3]{ba:le}.  Any simple
object in $\FF(U_q(\lie{sl}_2))$ is isomorphic to $ V_j$ for some $j$
satisfying
$$j = 
 \begin{cases} 
  0,1/2,1,\ldots            &\text{if $q$ is not a root of unity}\ ;
 \\
  0,1/2,1, \ldots, (r-2)/2    &\text{if $q = \exp(\frac{\pi i}{r})$}\ .
\end{cases}
$$
If $q = \exp( \pi i/r)$ the branching rule for the product is
$$ V_{j_1} \otimes V_{j_2} = \bigoplus_j V_j $$
where $j$ satisfies the quantum Clebsch-Gordan rules \eqref{qcg}.  If
$q$ is not a root of unity, the branching rule is that for
representations of $\lie{sl}_2$.  

\subsection{The $6j$ symbols}

$6j$ symbols for $\lie{sl}_2$ were introduced as a tool in atomic
spectroscopy by Racah \cite{ra:6j}, and then studied mathematically by
Wigner; see \cite{wi:gr}.  $6j$ symbols for $U_q(\lie{sl}_2)$ were
introduced by Kirillov and Reshetikhin \cite{ki:re}, who used them to
generalize the Jones knot invariant.  The material below can be found
in the book by Carter, Flath, and Saito \cite{ca:6j}.  If
$\Hom_{U_q(\sl_2)}(V_j,V_{j_1} \otimes V_{j_2})$ is non-trivial we fix
as a generator the map $\Y{j}{j_1}{j_2}$ defined in
\cite[p.96]{ca:6j}.  (A small difference is that in \cite{ca:6j} the
map is from $V_{\hh}^{\otimes j}$ to $ V_{\hh}^{\otimes j_1} \otimes
V_\hh^{\otimes j_2}$.)  If $\Hom_{U_q(\sl_2)}(V_j,V_{j_1} \otimes
V_{j_2})$ is trivial, define $\Y{j}{j_1}{j_2} = 0$.  Let $j_{ab}, 1
\leq a < b \leq 4$ be non-negative half-integers.  Define
$$ \YYL{j_{14}}{j_{13}}{j_{12}}{j_{23}}{j_{34}} =
\left(\Y{j_{13}}{j_{12}}{j_{23}} \otimes \id \right) \circ \Y{j_{14}}{j_{13}}{j_{34}} $$
and
$$ \YYR{j_{14}}{j_{24}}{j_{12}}{j_{23}}{j_{34}} = 
\left(\id \otimes \Y{j_{24}}{j_{23}}{j_{34}} \right) \circ \Y{j_{14}}{j_{12}}{j_{24}} $$
where $\id$ denotes the identity.  By associativity of the tensor
product we have a basis for $\Hom_{U_q(\sl_2)}(V_{j_{14}}, V_{j_{12}}
\otimes V_{j_{23}} \otimes V_{j_{34}})$ for each way of
parenthesizing:
$$ \Hom_{U_q(\sl_2)}(V_{j_{14}}, (V_{j_{12}} \otimes V_{j_{23}}) \otimes
V_{j_{34}}) = \bigoplus_{j_{13}} \C \YYL{j_{14}}{j_{13}}{j_{12}}{j_{23}}{j_{34}} $$
and
$$ \Hom_{U_q(\sl_2)}(V_{j_14},
V_{j_{12}} \otimes (V_{j_{23}} \otimes V_{j_{34}}))
=\bigoplus_{j_{24}} \C \YYR{j_{14}}{j_{24}}{j_{12}}{j_{23}}{j_{34}} $$
where the sum is over $j_{13}$ resp. $j_{24}$ satisfying the quantum
Clebsch-Gordan inequalities for each vertex.  The $6j$ symbol $\bsj
j_{12} & j_{23} & j_{13} \\ j_{34} & j_{14} & j_{24} \esj_{0} $ is
defined by
$$  \YYR{j_{14}}{j_{24}}{j_{12}}{j_{23}}{j_{34}} = 
\sum_{j_{13}} 
\bsj j_{12} & j_{23} & j_{13} \\ j_{34} & j_{14} & j_{24} \esj_{0} 
\YYL{j_{14}}{j_{13}}{j_{12}}{j_{23}}{j_{34}} .$$

Assume that $q$ is real or $|q| = 1$, so that the quantum integers are
real.  One can modify the definition slightly so that the symbols have
tetrahedral symmetry:
$$ 
\bsj j_{12} & j_{23} & j_{13} \\ j_{34} & j_{14} & j_{24} \esj 
= \frac{(-1)^{j_{12} + j_{23} + j_{34} + j_{14}}}{[2 j_{13} + 1]}
\sqrt{ \left|
 \frac{\theta(123) \theta(134)}
{ \theta(234) \theta (124)} \right| }
\bsj j_{12} & j_{23} & j_{13} \\ j_{34} & j_{14} & j_{24} \esj_{0} 
$$
where
$$ \theta(abc) =
\frac{ [j_{ab} + j_{bc} - j_{ac}]![j_{ab} - j_{bc} + j_{ac}]!
[-j_{ab} + j_{bc} + j_{ac}]![j_{ab} + j_{bc} + j_{ac} + 1]!}
{ (-1)^{j_{ab} + j_{bc} + j_{ac}} [2j_{ab}]![2j_{bc}]![2j_{ac}]!}$$
and the quantum factorial is defined by 
\begin{equation} \label{qfac} [n]! = [n][n-1]\ldots [1] .\end{equation}
That is, 
$$
\bsj j_{12} & j_{23} & j_{13} \\ j_{34} & j_{14} & j_{24} \esj
= 
\bsj j_{\sigma(1)\sigma(2)} & j_{\sigma(2)\sigma(3)} &
j_{\sigma(1)\sigma(3)}\\
 j_{\sigma(3)\sigma(4)} & j_{\sigma(1)\sigma(4)} & j_{\sigma(2)\sigma(4)} \esj
$$
for any $\sigma \in S_4$.  Equivalently, the $6j$ symbol is invariant
under the exchange of any two columns, and any two entries in the
first row with the corresponding entries in the second.  A word on
sign conventions: In the case $q$ arbitrary one defines the $6j$
symbols without taking the absolute value of the $\theta$ values.
This leads to a difference of sign $(-1)^{\sum_{a < b} j_{ab}}$ from
our convention, which was chosen because it agrees with the usual
definition of classical $6j$ symbols.  For readers familiar with spin
networks, the $6j$ symbol is the value of a tetrahedral graph with
labels $j_{ab}$, divided by the square root of the absolute value of
the product of the theta graphs for the faces.  This graph is {\em
dual} to the edge graph of the geometric tetrahedron, that is, edges
$j_{ab},j_{bc}, j_{ac}$ are joined at a vertex, for any $1 \leq a < b
< c \leq 4$.

The $6j$ symbols satisfy the {\em pentagon} or
{\em Biedenharn-Elliot relation}
\begin{equation} \label{pentagon}
\left\{
\begin{array}{ccc} j_{23} & j_{34} & j_{24} \\ 
                   j_{14} & j_{12} & j_{13} 
\end{array} \right\}
\left\{
\begin{array}{ccc} j_{23} & j_{34} & j_{24} \\ 
                   j_{45} & j_{25} & j_{35} 
\end{array} \right\} =
\end{equation}
\begin{equation*} \label{be2}
\sum_{j_{15}}(-1)^{z}[2j_{15}+1]
\left\{
\begin{array}{ccc} j_{13} & j_{34} & j_{14} \\ 
                   j_{45} & j_{15} & j_{35} 
\end{array} \right\}
\left\{
\begin{array}{ccc} j_{12} & j_{24} & j_{14} \\ 
                   j_{45} & j_{15} & j_{25} 
\end{array} \right\}
\left\{
\begin{array}{ccc} j_{12} & j_{23} & j_{13} \\ 
                   j_{35} & j_{15} & j_{25} 
\end{array} \right\}
\end{equation*}
where 
$ z=j_{12}+j_{13}+j_{14}+j_{15}+j_{23}+j_{24}+j_{25}+j_{34}+j_{35}+j_{45} .$
This identity is obtained by applying the definition of the $6j$
symbols to Figure 4.
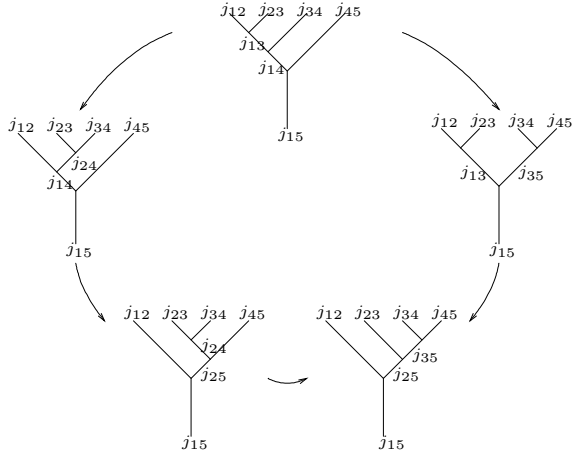
\begin{figure}[h]
\setlength{\unitlength}{0.00033333in}
\begingroup\makeatletter\ifx\SetFigFont\undefined%
\gdef\SetFigFont#1#2#3#4#5{%
  \reset@font\fontsize{#1}{#2pt}%
  \fontfamily{#3}\fontseries{#4}\fontshape{#5}%
  \selectfont}%
\fi\endgroup%
{\renewcommand{\dashlinestretch}{30}
\begin{picture}(8676,6963)(0,-10)
\path(3450,6804)(4350,5904)(4350,5004)
\path(4050,6804)(3750,6504)
\path(4650,6804)(4050,6204)
\path(4350,5904)(5250,6804)
\path(6750,2004)(5850,1104)(4950,2004)
\path(5850,1104)(5850,204)
\path(5550,2004)(6150,1404)
\path(6150,2004)(6450,1704)
\path(2550,2004)(3150,1404)
\path(3150,2004)(2850,1704)
\path(1950,2004)(2850,1104)(3750,2004)
\path(2850,1104)(2850,204)
\path(1950,4929)(1050,4029)(150,4929)
\path(1350,4929)(750,4329)
\path(750,4929)(1050,4629)
\path(1050,4029)(1050,3204)
\path(6750,5004)(7650,4104)(8550,5004)
\path(7350,5004)(7050,4704)
\path(7950,5004)(8250,4704)
\path(7650,4104)(7650,3204)
{\tiny 
\put(3300,6804){\makebox(0,0)[lb]{\smash{{{$j_{12}$}}}}}
\put(3900,6804){\makebox(0,0)[lb]{\smash{{{$j_{23}$}}}}}
\put(4500,6804){\makebox(0,0)[lb]{\smash{{{$j_{34}$}}}}}
\put(5100,6804){\makebox(0,0)[lb]{\smash{{{$j_{45}$}}}}}
\put(0,5004){\makebox(0,0)[lb]{\smash{{{$j_{12}$}}}}}
\put(600,5004){\makebox(0,0)[lb]{\smash{{{$j_{23}$}}}}}
\put(1200,5004){\makebox(0,0)[lb]{\smash{{{$j_{34}$}}}}}
\put(1800,5004){\makebox(0,0)[lb]{\smash{{{$j_{45}$}}}}}
\put(6600,5079){\makebox(0,0)[lb]{\smash{{{$j_{12}$}}}}}
\put(7200,5079){\makebox(0,0)[lb]{\smash{{{$j_{23}$}}}}}
\put(7800,5079){\makebox(0,0)[lb]{\smash{{{$j_{34}$}}}}}
\put(8400,5079){\makebox(0,0)[lb]{\smash{{{$j_{45}$}}}}}
\put(4800,2079){\makebox(0,0)[lb]{\smash{{{$j_{12}$}}}}}
\put(5400,2079){\makebox(0,0)[lb]{\smash{{{$j_{23}$}}}}}
\put(6000,2079){\makebox(0,0)[lb]{\smash{{{$j_{34}$}}}}}
\put(6600,2079){\makebox(0,0)[lb]{\smash{{{$j_{45}$}}}}}
\put(1800,2079){\makebox(0,0)[lb]{\smash{{{$j_{12}$}}}}}
\put(2400,2079){\makebox(0,0)[lb]{\smash{{{$j_{23}$}}}}}
\put(3000,2079){\makebox(0,0)[lb]{\smash{{{$j_{34}$}}}}}
\put(3600,2079){\makebox(0,0)[lb]{\smash{{{$j_{45}$}}}}}}
\path(1050,2904)(1050,2902)(1051,2899)
	(1052,2893)(1054,2883)(1056,2870)
	(1059,2854)(1063,2836)(1068,2815)
	(1073,2791)(1080,2767)(1087,2741)
	(1095,2713)(1104,2683)(1115,2652)
	(1127,2618)(1142,2581)(1159,2541)
	(1178,2499)(1200,2454)(1221,2413)
	(1242,2374)(1263,2337)(1283,2303)
	(1303,2271)(1322,2242)(1341,2214)
	(1359,2187)(1377,2163)(1394,2139)
	(1411,2116)(1427,2095)(1442,2076)
	(1456,2059)(1468,2043)(1478,2030)
	(1487,2020)(1500,2004)
\path(1401.045,2078.216)(1500.000,2004.000)(1447.612,2116.051)
\path(2550,6504)(2549,6504)(2547,6503)
	(2542,6501)(2535,6498)(2526,6494)
	(2513,6489)(2498,6483)(2480,6476)
	(2459,6467)(2436,6457)(2411,6445)
	(2384,6433)(2355,6420)(2326,6405)
	(2294,6389)(2262,6373)(2229,6354)
	(2194,6335)(2158,6313)(2120,6290)
	(2081,6265)(2039,6237)(1995,6206)
	(1949,6172)(1901,6135)(1851,6096)
	(1800,6054)(1753,6014)(1707,5973)
	(1663,5932)(1622,5893)(1582,5854)
	(1545,5817)(1510,5780)(1477,5745)
	(1446,5711)(1417,5677)(1388,5644)
	(1361,5612)(1335,5580)(1310,5550)
	(1287,5520)(1264,5491)(1243,5463)
	(1223,5437)(1205,5412)(1188,5390)
	(1173,5370)(1160,5353)(1150,5338)
	(1141,5327)(1135,5318)(1125,5304)
\path(1170.337,5419.085)(1125.000,5304.000)(1219.161,5384.211)
\path(7650,2904)(7650,2902)(7649,2899)
	(7648,2893)(7646,2883)(7644,2870)
	(7641,2854)(7637,2836)(7632,2815)
	(7627,2791)(7620,2767)(7613,2741)
	(7605,2713)(7596,2683)(7585,2652)
	(7573,2618)(7558,2581)(7541,2541)
	(7522,2499)(7500,2454)(7479,2413)
	(7458,2374)(7437,2337)(7417,2303)
	(7397,2271)(7378,2242)(7359,2214)
	(7341,2187)(7323,2163)(7306,2139)
	(7289,2116)(7273,2095)(7258,2076)
	(7244,2059)(7232,2043)(7222,2030)
	(7213,2020)(7200,2004)
\path(7252.388,2116.051)(7200.000,2004.000)(7298.955,2078.216)
\path(6150,6504)(6151,6504)(6153,6503)
	(6157,6501)(6164,6498)(6173,6494)
	(6185,6489)(6200,6483)(6217,6476)
	(6237,6467)(6259,6457)(6284,6445)
	(6310,6433)(6338,6420)(6367,6405)
	(6397,6389)(6429,6373)(6462,6354)
	(6496,6335)(6532,6313)(6570,6290)
	(6610,6265)(6653,6237)(6697,6206)
	(6745,6172)(6795,6135)(6847,6096)
	(6900,6054)(6947,6016)(6993,5978)
	(7037,5940)(7079,5902)(7120,5866)
	(7158,5831)(7195,5796)(7230,5762)
	(7263,5730)(7295,5698)(7325,5667)
	(7355,5636)(7383,5606)(7411,5576)
	(7437,5548)(7463,5520)(7487,5492)
	(7510,5466)(7532,5442)(7552,5418)
	(7571,5397)(7588,5377)(7602,5360)
	(7615,5345)(7626,5333)(7634,5323)
	(7640,5315)(7650,5304)
\path(7547.081,5372.613)(7650.000,5304.000)(7591.477,5412.973)
\path(4554.799,1025.027)(4650.000,1104.000)(4528.835,1079.118)
\path(4650,1104)(4625,1092)(4611,1086)
	(4594,1079)(4575,1071)(4554,1064)
	(4532,1056)(4508,1049)(4482,1043)
	(4453,1037)(4422,1033)(4387,1030)
	(4350,1029)(4313,1030)(4278,1033)
	(4247,1037)(4218,1043)(4192,1049)
	(4168,1056)(4146,1064)(4125,1071)
	(4106,1079)(4089,1086)(4075,1092)
	(4064,1097)(4056,1101)(4052,1103)(4050,1104)
{\tiny
\put(4200,4854){\makebox(0,0)[lb]{\smash{{{$j_{15}$}}}}}
\put(900,3054){\makebox(0,0)[lb]{\smash{{{$j_{15}$}}}}}
\put(2700,54){\makebox(0,0)[lb]{\smash{{{$j_{15}$}}}}}
\put(5775,54){\makebox(0,0)[lb]{\smash{{{$j_{15}$}}}}}
\put(7500,3054){\makebox(0,0)[lb]{\smash{{{$j_{15}$}}}}}
\put(3900,5904){\makebox(0,0)[lb]{\smash{{{$j_{14}$}}}}}
\put(3600,6279){\makebox(0,0)[lb]{\smash{{{$j_{13}$}}}}}
\put(7050,4254){\makebox(0,0)[lb]{\smash{{{$j_{13}$}}}}}
\put(7950,4254){\makebox(0,0)[lb]{\smash{{{$j_{35}$}}}}}
\put(6000,1104){\makebox(0,0)[lb]{\smash{{{$j_{25}$}}}}}
\put(6300,1404){\makebox(0,0)[lb]{\smash{{{$j_{35}$}}}}}
\put(3000,1104){\makebox(0,0)[lb]{\smash{{{$j_{25}$}}}}}
\put(3000,1554){\makebox(0,0)[lb]{\smash{{{$j_{24}$}}}}}
\put(600,4104){\makebox(0,0)[lb]{\smash{{{$j_{14}$}}}}}
\put(975,4404){\makebox(0,0)[lb]{\smash{{{$j_{24}$}}}}}}
\end{picture}
}
\caption{The pentagon identity}
\end{figure}

Racah gave a hypergeometric formula whose generalization to arbitrary
$q$ is
\begin{multline} \label{racah}
\bsj j_{12} & j_{23} & j_{13} \\ j_{34} & j_{14} & j_{24} \esj  =
\Delta(123) \Delta(134) \Delta(234) \Delta (124) \sum (-1)^z
[z+1]! f(z)^{-1} \end{multline}
where
\begin{multline}
f(z) = [z - j_{12} -j_{23} - j_{13}]!  [z - j_{13} - j_{34} - j_{14}]!
[z - j_{23} - j_{34} - j_{24}]! 
\\ [z - j_{12} - j_{24} - j_{14}]!
[j_{12} + j_{23} + j_{34} + j_{14} - z]!  [j_{12} + j_{13} + j_{34} +
j_{24} - z]!  
\\ [j_{23} + j_{13} + j_{14} + j_{24} - z]!,
\end{multline}
$$ \Delta(abc)
= \left( \frac{[j_{ab} + j_{bc} - j_{ac}]![j_{ab} - j_{bc} + j_{ac}]!
[-j_{ab} + j_{bc} + j_{ac}]!}{[j_{ab} + j_{bc} + j_{ac} + 1]!}
\right)^{1/2} $$
and the sum is over integers $z$ such that the factorials in $f(z)$
are defined.  Because a large amount of cancellation occurs,
\eqref{racah} is useful for studying asymptotics only in special
cases.

\subsection{Asymptotics for $q = 1$}

The asymptotics of the $6j$ symbols for $q = 1$ were considered by
Wigner \cite{wi:gr} as follows.  Consider a system of particles A,B,C
with total angular momenta $j_{12},j_{23}, j_{34}$.  Given that the
total angular momentum of the combined AB system is $j_{13}$, and the
total angular momentum of the combined ABC system is $j_{14}$, the
probability of measuring $j_{24}$ for the total angular momentum of
the combined BC system is the square of the $6j$ symbol.  Wigner
computed the classical probability and arrived at the conjecture
$$ \bsj j_{12} & j_{23} & j_{13} \\ j_{34} & j_{14} &
j_{24} \esj^2 \approx \frac{1}{24 \pi \Vol(\tau)} $$
where $\tau$ is the Euclidean tetrahedron with lengths $l_{ab} =
j_{ab} + \hh$, if it exists, and the approximation is meant to hold
after averaging out local oscillations.  The Wigner conjecture was
refined by Ponzano and Regge \cite{po:6j} as follows.  Let
$\theta_{ab}$ denote the exterior dihedral angles of $\tau$, and
$$ \phi = \sum_{a < b} l_{ab} \theta_{ab} .$$
Ponzano and Regge conjectured in the semiclassical limit
\begin{equation} \label{pr}  
\bsj j_{12} & j_{23} & j_{13} \\    j_{34} & j_{14} & j_{24} \esj 
\sim \frac{\cos(\phi + \pi/4)}{\sqrt{12 \pi \Vol(\tau)}} \end{equation}
assuming that $\tau$ exists and is non-degenerate, and also formulas
for degenerate tetrahedra with non-degenerate faces and the
classically forbidden case.  The Ponzano-Regge conjectures were
further refined by Schulten and Gordon \cite{sch:semi}, see Section 4
below.  The mathematical meaning of the semiclassical limit is not
entirely clear.  An obvious interpretation (but perhaps not the only
one) is the asymptotics of the sequence $ \bsj kj_{12} & kj_{23} &
kj_{13} \\ kj_{34} & kj_{14} & kj_{24} \esj, k \to \infty $.  In this
setting the Ponzano-Regge formula \eqref{pr} has recently been proved
by J. Roberts \cite{ro:qu}, see also \cite{ch:as},\cite{wi:as}.  The
conjectures for the other cases follow from the recursion in
\cite{sch:semi} and its mathematical justification by Geronimo, Bruno
and Assche \cite[Theorem 3.8]{ge:tp}.  The remaining cases (c)-(f) can
be handled using Stirling's formula applied to \eqref{racah}, see also
Sections 6 and 7 below.

\subsection{Non-Euclidean tetrahedra}

Let $E^n,S^n, H^n$ denote Euclidean, spherical, or hyperbolic
$n$-space respectively.  Let $\sigma$ be an $n$-simplex in $E^n,S^n$
or $H^n$.  Let $l_{ab}, \ 0 \leq a,b \leq n$ denote the lengths of the
edges, and $l$ the matrix with coefficients $l_{ab}$.  The study of
which edge lengths occur is called {\em distance geometry} in
Blumenthal \cite{bl:di}.  Let $ G_0(l)$ (also denoted $G_0(\sigma)$
denote the {\em Cayley-Menger} $(n+2) \times (n+2)$-matrix (actually a
slight modification) written in block-diagonal form
\begin{equation} \label{Cayley}
G_0(l) = \left[ \begin{array}{rrrrrrr} 0 & 1 & 1 & 1 & \ldots & 1 \\ 
1 & 0 & -\hh l_{01}^2 & - \hh l_{02}^2 & \ldots & - \hh l_{0n}^2 \\
1 & - \hh l_{01}^2  & 0 & - \hh l_{12}^2 & \ldots & - \hh l_{1n}^2 \\
\vdots & \vdots & \vdots & & & \vdots \\
1 & -\hh l_{0n}^2 & - \hh l_{1n}^2 & \ldots & & 0  
\end{array} \right].
\end{equation} 
For a spherical or hyperbolic simplex $\sigma$ we define the {\em Gram
matrices} $G_\pm(l)$ (also denoted $G_\pm(\sigma)$) by
$$G_+(l)_{ab} = \cos(l_{ab}), \ \ G_-(l)_{ab} = - \cosh(l_{ab}) $$
and $|G_\pm(l)|$ its determinant.  We will need later the following
facts:

\begin{proposition} \label{nonEuclid} 
\begin{enumerate}
\item If an $n$-simplex $\sigma$ with edge lengths $l_{ab}$ exists in
$E^n$, then $\sigma$ has volume given by $ (n!\Vol(\sigma))^2 = - |G_0(l)|  .$
\item An $n$-simplex $\sigma$ with edge lengths $l_{ab}$ exists in
$E^n$ if and only if all principal minors in $G_0(l)$ containing the
$11$-entry are non-positive.  In this case, $\sigma$ is non-degenerate
if and only if $|G_0(l)| < 0$.
\item An $n$-simplex with edge lengths $l_{ab}$ exists in $S^n$ if and
only if $G_+(l)$ is positive semidefinite.  In this case, the simplex
is non-degenerate if and only if $|G_+(l)| > 0$.
\item An $n$-simplex with edge lengths $l_{ab}$ exists in $H^n$ if
and only if all principal minors of $G_-(l)$ are non-positive.  In
this case, the simplex is non-degenerate if and only if
$|G_-(l)| < 0$.
\item If a tetrahedron $\tau$ with edge lengths $l_{ab}$ exists, it is
unique up to isometry (not necessarily orientation-preserving) in
$E^n$ resp. $S^n, H^n$.  If $\tau$ is non-degenerate, there are two
equivalence classes of tetrahedra up to orientation-preserving
isometry, called {\em mirror tetrahedra}.
 \item \label{Gram2} In the case $n = 1$ we have
$$ |G_0(l)| = -l_{12}^2, \ \ \ |G_+(l)| =
\sin^2(l_{12}), \ \ \ |G_-(l)| = -\sinh^2(l_{12}) .$$
\item \label{Gram3} In the case $n = 2$, there are factorizations,
$$ |G_0(l)| = -4 s(s-l_{12})(s-l_{23})(s - l_{13}) $$
$$ |G_+(l)| = 4 \sin (s) \sin(s-l_{12})
\sin(s-l_{23}) \sin (s - l_{13}) $$
$$ |G_-(l)| = -4
 \sinh (s) \sinh(s-l_{12})
\sinh(s-l_{23}) \sinh (s - l_{13}) $$
where $s = (l_{12} + l_{23} + l_{13})/2$ is the semiperimeter.  (These
quantities appear in the numerators of the Heron area formulas.)
\item 
A Euclidean, resp. spherical, resp. hyperbolic triangle exists if and
only if \eqref{tri}, resp. \eqref{triq} and \eqref{tri},
resp. \eqref{tri} hold:
\begin{equation} \label{tri} 
l_{ac} \leq l_{ab} + l_{bc}, \ \ l_{ab} \leq l_{ac} + l_{bc}, \ \
l_{bc} \leq l_{ab} + l_{ac} .\end{equation}
\begin{equation} \label{triq} 
l_{ab} + l_{bc} + l_{ac} \leq 2\pi
.\end{equation}
\item In the case $n = 3$, a tetrahedron $\tau$ with edge lengths
$l_{ab}$ exists in $E^3$ resp. $S^3,H^3$ if and only if $l_{ab}$
satisfy the triangle inequalities for each face and
$$\det G_0(l) \leq 0, \ \text{resp.}  \det G_+(l) \geq 0, \
 \ \text{resp.} \det G_-(l) \leq 0 .$$
\item In the limit of small lengths, the determinants of $G_\pm $
approach the determinant of the Cayley-Menger matrix, that is,
\label{GG0}
$$ \lim_{\eps \to 0} \frac{\det G_\pm(\eps l)}{ \eps^{2n}} = \mp
 \det G_0(l).
$$
%
\item \label{equal}
A spherical, resp. hyperbolic $n$-simplex with edge lengths $l$
exists if and only if $0 \leq l \leq \arccos(-1/n)$, resp $0 \leq l$.
\item The volumes of non-degenerate $n$-simplices $\sigma$ in
$E^n, S^n$ or $H^n$ satisfy an identity due to Schl\"afli
$$ (n-1) \kappa \d \Vol(\sigma) = \sum \Area(F) \d \theta_F $$
where the sum is over $(n-2)$-dimensional faces $F$ of the simplex
$\sigma$, $\theta_F$ is the exterior dihedral angle, and $\kappa =
0,1,-1$ resp. is the curvature.  
\item \label{schlafli} 
For non-degenerate tetrahedra $\tau$ in $S^3,H^3,or E^3$,
$$ 2 \kappa \d \Vol(\tau(s)) = \sum_{a < b} l_{ab} \d \theta_{ab}.$$
\item The derivative of dihedral angle with respect to opposite
edge length is 
\label{deriv} 
$$  \left( \frac{\partial
\theta_{ab}}{\partial l_{cd}} \right)^{-1}
= \left \{ \begin{array}{ll}
- ||G_0(l)||^{-1/2}{l_{cd} l_{ab}} & \text{in
  the 
Euclidean case} \\
- |G_+(l)|^{-1/2}{\sin(l_{cd}) \sin(l_{ab})}
& \text{in the spherical case} \\
- ||G_-(l)||^{-1/2}{\sinh(l_{cd})
  \sinh(l_{ab})}
& \text{in the hyperbolic case} 
\end{array} \right\} $$
\end{enumerate}
\end{proposition}

\begin{proof} (a-d) are all discussed in \cite{bl:di}.  Note that the
non-degenerate cases of (c,d) can be treated uniformly using a theorem
of Jacobi \cite[p.303]{gant:ma}, as was pointed out to us by F. Luo.
The case of degenerate hyperbolic simplices is slightly problematic,
and we do not know a reference; it can be proved by induction on $n$.
We will not use this case in the paper.  (e-j) are all
straight-forward.  (k) follows from the fact that for all edge lengths
equal to $l$ we have 
$$ | G_+(l) | = (1-\cos(l))^{n-1}(1+n\cos(l)), 
\ \ |G_-(l)| = - (\cosh(l) - 1)^{n-1}(n\cosh(l) + 1) .$$  
For (l),(m) see \cite[p. 281]{mil:col}.  We remark that a formula for
the volume of a hyperbolic tetrahedron in terms of its edge lengths is
given in Murakami-Ushijima \cite{mu:vo}, see also Mohanty
\cite{mo:th}.  (n) $|| G_0(l) ||$ refers to the absolute value of the
determinant.  The Euclidean case is a computation of Wigner.  We prove
the claim for the spherical case.  Let $V_1,\ldots,V_4$ be matrices in
$SU(2) \cong S^3$, and $E_{ij} = V_i V_j^{-1}$.  We may
assume 
$$V_a = \bma{rr} 1 & 0 \\ 0 & 1 \ema, \ \ V_b = 
\bma{rr} e^{il_{ab}} & 0 \\ 0 & e^{-il_{ab}} \ema .$$
Let $D_\theta = \diag(\exp(i \theta/2),\exp(- i\theta/2))$.  Rotation
around $E_{ab}$ by angle $\theta$ corresponds to conjugation by
$D_\theta$.  We compute
$$ \frac{\d}{\d \theta} \vert_{\theta = 0 } \Tr(D_\theta V_c
 D_\theta^{-1} V_d^{-1}) = -2 \det \left[ \begin{array}{cc} V_c^3 & V_d^3 \\
 V_c^4 & V_d^4\end{array} \right]
$$
where $V_i^j, i,j = 1,2,3,4$ are the components of $V_i$ as a vector
in $\R^4$, using the diffeomorphism $SU(2) \to S^3$.  Let $V$
be the matrix with coefficients $V_i^j$.  Since 
$$ (V^T V)_{ij} = V_i \cdot V_j = \cos(l_{ij}) = G_+(l)_{ij} $$
we have $\det(V)^2 = \det(G_+(l))$.  Since $V$ is in block-diagonal form, 
\begin{eqnarray}
\det(V) &=& \det \left[ \begin{array}{cc} V_a^1 & V_b^1 \\ V_a^2 & V_b^2
\end{array} \right] 
 \det \left[ \begin{array}{cc} V_c^3 & V_d^3 \\ V_c^4 & V_d^4
\end{array} \right]   \\
&=& - \left[ \begin{array}{cc} 1 & \cos(l_{ab}) \\ 0 & \sin(l_{ab})
   \end{array} \right]   \hh \frac{\d \Tr(E_{cd}))}{\d \theta} \\
&=& - \sin(l_{ab}) \sin(l_{cd}) \left( \frac{\partial \theta_{ab}}{\partial
l_{cd}} \right)^{-1}
\end{eqnarray}
which proves the spherical part of the lemma.  The hyperbolic case is
similar, using singular values instead of eigenvalues.
\end{proof}

\vskip .1in

The following gives a particularly useful formula for the dihedral
angle of a spherical tetrahedron. Let $v_0, \ldots, v_n$ in
$\mathbb{R}^n$.  For any $\{ i_1,\ldots,i_k \} \subset \{ 1,\ldots,
n\}$ we denote by $\sigma(i_1,i_2,\ldots,i_k)$ the $n-k$ simplex
spanned by $v_j$ for $j \neq i_1,\ldots,i_k$.
\begin{lemma} \label{nsimplex} For $a,b \in \{1,\ldots,n \}$ distinct,
let $\theta_{ab}$ be the exterior dihedral angle between
$\sigma(a),\sigma(b)$.
$ n \Vol (\sigma) \Vol (\sigma(a,b)) = (n-1) \Vol (\sigma(a))
\Vol (\sigma(b)) \sin(\theta_{ab}) .$
\end{lemma}

\begin{proof}
Let $h_{a}$ denote the orthogonal projection of $v_a$ to onto the
orthogonal complement of the span of $\sigma(a,b)$ for $i=a,b$.
$$n (n-1) \Vol (\sigma)= \Vert h_{a} \Vert \Vert h_b \Vert \sin \theta
\Vol(\sigma(a,b)).$$
Using
$(n-1) \Vol (\sigma(a))= \Vert h_{b} \Vert \Vol(\sigma(a,b)) $
and similarly for $\sigma(b)$ proves the claim.
\end{proof}

\begin{corollary} \label{angle}
Let $\tau$ be a spherical tetrahedron with vertices $v_a,v_b,v_c,v_d
\in S^3$.
\begin{equation}\label{vol}
| G_+(\tau)|^{\hh} \sin l_{ab}=|G_+(\tau(d))|^{\hh} |G_+(\tau(c))|^{\hh} \sin
\theta_{ab},
\end{equation}
where $G_+(\tau)$ denotes the Gram matrix with vertices $\tau$ etc.
Furthermore,
\begin{equation}\label{dgram}
2 |G_+(\tau(d))|^{1/2} |G_+(\tau(c))|^{1/2}\cos \theta_{ab} = 
\frac{d |G_+(\tau)|}{d \cos l_{cd}}.
\end{equation}
\end{corollary}

\begin{proof}
\eqref{vol} is a special case of Lemma \ref{nsimplex}.  \eqref{dgram}
is the derivative of \eqref{vol} with respect to $\cos l_{cd}$.  We
start with squaring \eqref{vol}:
\begin{equation}\label{squared}
|G_+(\tau)|\sin^2 l_{ab}=|G_+(\tau(d))||G_+(\tau(c))|\sin^2 \theta_{ab}.
\end{equation}
Taking the derivative of the left hand side gives 
$$\frac{d |G_+(\tau)|}{\d \cos l_{cd}}\sin^2
l_{ab}.$$
The derivative of the right hand side is
$$|G_+(\tau(d))||G_+(\tau(c))|2\sin \theta_{ab} \cos \theta_{ab} \frac{\d \theta_{ab}}{\d \cos
  l_{cd}}.$$
Substituting \eqref{vol} in \eqref{squared}, one obtains 
\begin{equation}\label{rhs}
|G_+(\tau(d))|^{1/2}|G_+(\tau(c))|^{1/2}2\cos \theta_{ab}
|G_+(\tau)|^{1/2}\sin l_{ab}
\frac{\d\theta_{ab}}{\d\cos l_{cd}}.
\end{equation}

We compute 

$$\frac{\d\theta_{ab}}{\d \cos
  l_{cd}}=\frac{\d\theta_{ab}}{\d l_{cd}}\frac{\d l_{cd}}{\d\cos l_{cd}}
= -\left(\frac{\sin l_{ab} \sin
  l_{cd}}{ |G_+(\tau)|^{1/2}}\right)\left(\frac{
-1}{\sin l_{cd}}\right).$$
Substitute this derivative in the equation \eqref{rhs}, we get
$$2\pi |G_+(\tau(d))|^{1/2} |G_+(\tau(c))|^{1/2} \cos \theta_{ab} \sin^2 l_{ab}.$$
Thus, the derivative of the whole equation is
$$\frac{\d |G_+(\tau)|}{\d \cos l_{cd}}=
2|G_+(\tau(d))|^{1/2}|G_+(\tau(c))|^{1/2}\cos \theta_{ab},$$
as required.  
\end{proof}

We remark that 
$$ \frac{\d |G_+(\tau)|}{\d \cos l_{cd}}=  -2 
\left|    \begin{array}{ccc} 1 & \cos(l_{ab}) & \cos(l_{ad}) \\
                   \cos(l_{ab}) &     1       & \cos(l_{bd}) \\
                   \cos(l_{ac}) & \cos(l_{bc}) & \cos(l_{cd}) \end{array}
\right|. 
$$
\section{Geometry of conformal blocks}

In this section we explain the geometry of the conjecture
\eqref{conj}.  Let $p_{1},p_{2},p_{3},p_{4}$ be distinct points on
$\P^1$.  Given $l_{a(a+1)} \in [0,\pi], a = 1,\ldots,4$ let 
$$ t_a = \bma{cc}  e^{i l_{a(a+1)}} & 0 \\  0 & e^{- i
l_{a(a+1)}} \ema .$$
Let $\M = \M(l_{12},l_{23},l_{34},l_{41})$ denote the moduli space of
flat $SU(2)$-bundles on $\P^1 - \{ p_{1},p_{2},p_{3},p_{4} \}$ with
holonomy around $p_a$ conjugate to $t_a$.  The fundamental group of
$\P^1 - \{ p_{1},p_{2},p_{3},p_{4} \}$ is generated by the loops
$\gamma_j$ around $p_j$, with the single relation $\gamma_1\gamma_2
\gamma_3 \gamma_4 = 1$.  Hence
$$ \M = \{ (e_1,e_2,e_3,e_4) \in SU(2)^4, \ e_1e_2e_3e_4 = I, \ e_a
\sim t_a,\ \ a = 1,2,3,4 \}/SU(2) .$$
Under the diffeomorphism $SU(2) \to S^3$, the set of matrices
conjugate to $t_a$ is mapped to the set of points in $S^3$ at distance
$l_{a(a+1)}$ from $(1,0,0)$.  The $4$-tuple
$$ v_1 = I, \ \ v_2 = e_1, \ \ v_3 = e_1e_2, \ \ v_4 = e_1e_2e_3 $$
maps to the vertices of a closed $4$-gon in $S^3$ with edge lengths
$l_{a(a+1)}$.  From now on, we assume $l_{a(a+1)} < \pi$.  Then $\M =
\M(l_{12},l_{23},l_{34},l_{41})$ is the moduli space of closed
$4$-gons in $S^3$ with edge lengths $l_{a(a+1)}$, where we identify
two $4$-gons if they are related by an isometry of $S^3$.  $\M$ is
non-empty if and only if the inequalities
$$l_{14} \leq l_{12} + l_{23} + l_{34}, \ \ \ l_{12} + l_{23} + l_{34}
\leq 2 \pi + l_{14} $$
and their cyclic permutations hold; this can be seen either directly
or from \cite{ag:ei},\cite{bl:ip} by computing Gromov-Witten
invariants of $\P^1$.  If one of these equalities holds with equality,
then the $4$-gon is contained in a geodesic on $S^2$ and $\M$ consists
of a single point.  If one of the lengths is zero then $\M$ is single
point.  Otherwise, $\M$ is homeomorphic to $S^2$.  If in addition $
\pm l_{12} \pm l_{23} \pm l_{34} \pm l_{14}$ is not a multiple of
$2\pi$ for all choices of sign then $\M$ is non-singular and
symplectomorphic to $S^2$.

There are two canonical coordinate systems on $\M$, described as
follows.  Let $\lambda_{13}$, resp. $\lambda_{24}$, denote the
distance from the vertex $v_{1}$ to $v_3$, resp. $v_2$ to $v_4$.  In
the language of flat bundles, the functions $\lambda_{13}$,
$\lambda_{24}$ measure the eigenvalues of the holonomy of the
connection around a circle separating the four points into two groups
of two.  Let $\theta_{13}$, resp.  $\theta_{24}$, denote the
corresponding dihedral angle.  Let $\M^*$ denote the subset of
$4$-gons that are not contained in a geodesic, that is, the subset of
irreducible bundles.  Let
$$\M_{13}^* \subset \M^*, \ \ \M_{24}^* \subset \M^* $$ 
denote the subset on which $0 < \lambda_{13} < \pi$ and the faces
containing the edge $v_1v_3$ are non-degenerate, resp. $0 <
\lambda_{24} < \pi$ and the faces containing $v_2v_4$ are
non-degenerate.  A special case of a result of Goldman \cite{go:in}
states that $\lambda_{13},\theta_{13}$, respectively
$\lambda_{24},\theta_{24}$, are smooth on $\M_{13}^*$, resp.
$\M_{24}^*$ and form a set of action-angle coordinates.  That is,
$$ \omega |_{\M_{13}^*} =  \d \theta_{13} \wedge \d \lambda_{13}, \ \
\omega |_{\M_{24}^*} = \d \theta_{24} \wedge \d \lambda_{24}.$$
If
$l_{12} = l_{23}$ and $l_{34} = l_{41}$ then $\lambda_{13}^{-1}(0)$ is
a subset of $\M$ homeomorphic to $S^1$, meeting the complement of
$\M^*$ in the two configurations where all vertices are colinear.  Let
$l_{13}^{\min}$ denote the minimum value of $\lambda_{13}$.  If
$l_{13}^{\min} \neq 0$, the set defined by $\lambda_{13} =
l_{13}^{\min}$ is a point.  If $l_{13}^{\min} = 0$, the set
$\lambda_{13} = l_{13}^{\min}$ is a circle.

By a theorem of Mehta-Seshadri \cite{ms:pb} and Furuta-Steer
\cite{fu:se}, $\M$ has the structure of a normal projective variety.
By a theorem of Pauly \cite{pa:co}, $\M$ has a positive line bundle
$\L$ whose sections $H^0(\M,\L^{r-2})$ may be identified with the
space of genus zero WZW conformal blocks at level $r-2$ with weights
$\mu_1,\mu_2,\mu_3,\mu_4$.  Let $\ol{\mathcal{C}}_{0,4}$ denote the
moduli space of stable $4$-pointed genus zero curves; the open subset
${\mathcal{C}}_{0,4}$ of $4$-pointed smooth curves is an open subset
identified with $\C - \{0,1\}$ by the map $ (\P^1;0,1,\infty,z)
\mapsto z$.  Let $X_{13}$, resp. $X_{24}$, denote the point in
$\ol{\mathcal{C}}_{0,4}$ represented by a nodal curve with two
components and points $p_{1},p_{2}$ on one component and $p_{3},p_{4}$
on the other, resp. $p_{1},p_{4}$ on one component and $p_{2},p_{3}$
on the other.
\begin{figure}[h]
\setlength{\unitlength}{0.00020333in}
\begingroup\makeatletter\ifx\SetFigFont\undefined%
\gdef\SetFigFont#1#2#3#4#5{%
  \reset@font\fontsize{#1}{#2pt}%
  \fontfamily{#3}\fontseries{#4}\fontshape{#5}%
  \selectfont}%
\fi\endgroup%
{\renewcommand{\dashlinestretch}{30}
\begin{picture}(16815,10171)(0,-10)
\put(7433,7987){\ellipse{4322}{4322}}
\put(6476,6792){\blacken\ellipse{120}{120}}
\put(6476,6792){\ellipse{120}{120}}
\put(8628,6792){\blacken\ellipse{120}{120}}
\put(8628,6792){\ellipse{120}{120}}
\put(6476,9183){\blacken\ellipse{120}{120}}
\put(6476,9183){\ellipse{120}{120}}
\put(8628,9183){\blacken\ellipse{120}{120}}
\put(8628,9183){\ellipse{120}{120}}
\put(6296,6493){\makebox(0,0)[lb]{\smash{{{\SetFigFont{10}{12.0}{\rmdefault}{\mddefault}{\updefault}$p_3$}}}}}
\put(8449,6493){\makebox(0,0)[lb]{\smash{{{\SetFigFont{10}{12.0}{\rmdefault}{\mddefault}{\updefault}$p_4$}}}}}
\put(6237,9363){\makebox(0,0)[lb]{\smash{{{\SetFigFont{10}{12.0}{\rmdefault}{\mddefault}{\updefault}$p_1$}}}}}
\put(8389,9363){\makebox(0,0)[lb]{\smash{{{\SetFigFont{10}{12.0}{\rmdefault}{\mddefault}{\updefault}$p_2$}}}}}
\put(13170,7800){\ellipse{2400}{2400}}
\put(15570,7800){\ellipse{2474}{2474}}
\texture{55888888 88555555 5522a222 a2555555 55888888 88555555 552a2a2a 2a555555 
	55888888 88555555 55a222a2 22555555 55888888 88555555 552a2a2a 2a555555 
	55888888 88555555 5522a222 a2555555 55888888 88555555 552a2a2a 2a555555 
	55888888 88555555 55a222a2 22555555 55888888 88555555 552a2a2a 2a555555 }
\put(13170,8400){\shade\ellipse{150}{150}}
\put(13170,8400){\ellipse{150}{150}}
\put(13170,7200){\shade\ellipse{150}{150}}
\put(13170,7200){\ellipse{150}{150}}
\put(15570,8400){\shade\ellipse{150}{150}}
\put(15570,8400){\ellipse{150}{150}}
\put(15570,7200){\shade\ellipse{150}{150}}
\put(15570,7200){\ellipse{150}{150}}
\put(12570,8400){\makebox(0,0)[lb]{\smash{{{\SetFigFont{12}{14.4}{\rmdefault}{\mddefault}{\updefault}$p_1$}}}}}
\put(12570,7200){\makebox(0,0)[lb]{\smash{{{\SetFigFont{12}{14.4}{\rmdefault}{\mddefault}{\updefault}$p_4$}}}}}
\put(15870,8400){\makebox(0,0)[lb]{\smash{{{\SetFigFont{12}{14.4}{\rmdefault}{\mddefault}{\updefault}$p_2$}}}}}
\put(15870,7200){\makebox(0,0)[lb]{\smash{{{\SetFigFont{12}{14.4}{\rmdefault}{\mddefault}{\updefault}$p_3$}}}}}
\put(1245,8850){\ellipse{2400}{2400}}
\put(1245,6450){\ellipse{2474}{2474}}
\put(645,8850){\shade\ellipse{150}{150}}
\put(645,8850){\ellipse{150}{150}}
\put(1845,8850){\shade\ellipse{150}{150}}
\put(1845,8850){\ellipse{150}{150}}
\put(645,6450){\shade\ellipse{150}{150}}
\put(645,6450){\ellipse{150}{150}}
\put(1845,6450){\shade\ellipse{150}{150}}
\put(1845,6450){\ellipse{150}{150}}
\path(4545,7800)(3045,7800)
\path(4545,7800)(3045,7800)
\path(3165.000,7830.000)(3045.000,7800.000)(3165.000,7770.000)
\path(10170,7800)(11670,7800)
\path(10170,7800)(11670,7800)
\path(11550.000,7770.000)(11670.000,7800.000)(11550.000,7830.000)
\put(14045,5575){\makebox(0,0)[lb]{\smash{{{\SetFigFont{12}{14.4}{\rmdefault}{\mddefault}{\updefault}$X_{24}$}}}}}
\put(495,9150){\makebox(0,0)[lb]{\smash{{{\SetFigFont{12}{14.4}{\rmdefault}{\mddefault}{\updefault}$p_1$}}}}}
\put(1695,9150){\makebox(0,0)[lb]{\smash{{{\SetFigFont{12}{14.4}{\rmdefault}{\mddefault}{\updefault}$p_2$}}}}}
\put(1695,6750){\makebox(0,0)[lb]{\smash{{{\SetFigFont{12}{14.4}{\rmdefault}{\mddefault}{\updefault}$p_3$}}}}}
\put(495,6750){\makebox(0,0)[lb]{\smash{{{\SetFigFont{12}{14.4}{\rmdefault}{\mddefault}{\updefault}$p_4$}}}}}
\put(570,4575){\makebox(0,0)[lb]{\smash{{{\SetFigFont{12}{14.4}{\rmdefault}{\mddefault}{\updefault}$X_{13}$}}}}}
\end{picture}
}
\vskip -1in
\caption{Two degenerations of
  $(\P^1;p_{1},p_{2},p_{3},p_{4})$ in $\ol{\mathcal{C}}_{0,4}$.}
\end{figure}
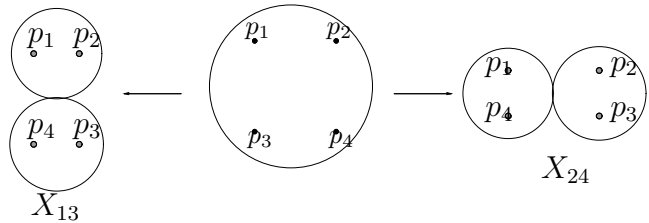
Let $\mathcal{V}$ denote the bundle with fiber
$
\mathcal{V}_{(\P^1,p_{1},p_{2},p_{3},p_{4})}
= H^0(\L^{r-2}(\P^1,p_1,p_2,p_3,p_4)) .$
According to the Khizhnik-Zamolodchikov, see
\cite{ts:cf},\cite{pa:co}, the bundle $\mathcal{V}$ has a projectively
flat connection $\nabla$.  Furthermore, the fiber over $X_{13}$,
resp. $X_{24}$ has a canonical projective basis $\{ s_{13}(j_{13}) \}
$ resp. $\{ s_{24}(j_{24}) \}$.  The $6j$ symbols are the entries in
the change of basis matrix for the bases $\{ s_{13}(j_{13}) \}, \{
s_{24}(j_{24}) \}$ after parallel transport to a common point in
$\ol{C}_{0,4}$.  (However, we only know how to give a canonical
normalization to the basis elements using the quantum group
definition.)  Define
$$ \Lambda_{13} = \lambda_{13}^{-1} \left(\frac{j_{13}}{r-2} \right), \ \
\Lambda_{24} = \lambda_{24}^{-1} \left(\frac{j_{24}}{r-2} \right) .$$
$\Lambda_{13},\Lambda_{24}$ are {\em Bohr-Sommerfeld} Lagrangians.
That is, the restriction of $\L^{r-2}$ to $\Lambda_{ab}$ has a flat
section $s^\infty_{ab}(j_{ab})$.  Extend the flat section to a section
$s^\infty_{ab}(j_{ab})$ defined in a neighborhood by parallel
transport in the perpendicular directions, followed by multiplication
by the Gaussian function $ \exp( - 4 \pi (r-2) d_{ab}^2) $ where
$d_{ab}$ is the geodesic distance to $\Lambda_{ab}$.
$s_\infty(j_{ab})$ is not in general holomorphic.  We conjecture that
as $k \to \infty$, after replacing $(r-2)$ with $k(r-2)$ one has 
\begin{equation} \label{scar}
s_{k,ab}(kj_{ab}) \sim s^\infty_{k,ab}(kj_{ab}) .\end{equation}
Part of the problem is to give a precise sense in which this
asymptotics holds.

We briefly discuss a possible approach to proving \eqref{scar}.  The
vector field $\frac{\partial}{\partial \theta_{13}}$ has (where
defined) a canonical lift to a connection-preserving vector field on
$\L$.  However, $\frac{\partial}{\partial \theta_{13}}$ does not
preserve the holomorphic structure, and so even if globally defined
does not define a circle action on $H^0(\L^{r-2})$.  However, we
conjecture that in the limit $(\P^1,p_{1},p_{2},p_{3},p_{4}) \to
X_{13}$ in $\ol{\mathcal{C}}_{0,4}$, the complex structure is
invariant; see Daskalopolous-Wentworth \cite{dask:fa1}.  This would
imply that the Verlinde sections $s_{13}(j_{13})$ are asymptotically
eigenvectors for the circle action with weight $j_{13}$.  \eqref{scar}
would follow from standard arguments discussed in \cite{ro:qu}.  We
also remark that the leading order part of $\nabla$ as discussed in
\cite[(4.23)]{ax:ge}) is a first order equation and should not change
the large $k$ asymptotics.

The determinant line bundle $\L$ is conjectured to admit a Hermitian
metric so that the induced Hermitian metric on $\mathcal{V}$ is
compatible with the projectively flat connection \cite{ax:ge}, see
also \cite{yosh:abel},\cite{ki:in}.  It would follow that the $6j$
symbol equals the Hermitian pairing of the sections $(s_{13}(j_{13}),
s_{24}(j_{24}))$.  Since $s_{k,ab}(kj_{ab})$ are asymptotically
Gaussian concentrating to $\Lambda_{ab}$, one naively expects the
pairing to have asymptotics
\begin{equation} \label{naive} 
\sum_{\tau \in \Lambda_{13} \cap \Lambda_{24}} \frac{(k(r-2))^{-1}
(s_{k,13}(kj_{13})(\tau), s_{k,24}(kj_{24})(\tau)) e^{\pm \pi i/4}}{ \sqrt{
(k\omega)_\tau( \ppthot,\ppthtf)}} .\end{equation}
We denote the two points of intersection $\tau_+,\tau_-$.  The
pairings $\omega_{\tau_\pm}( \ppthot,\ppthtf) $ are by \ref{nonEuclid}
\eqref{deriv}
$$ \omega_{\tau_\pm} \left( \ppthot,\ppthtf \right) = \frac{\partial
l_{13}}{\partial \theta_{24}} = - \frac{|G_+(l)|^{1/2}}{\sin(\mu_{cd})
\sin(\mu_{ab})}.
$$
Moving along the manifold $\Lambda_{13}$ from $\tau_+$ to $\tau_-$ and
then $\Lambda_{24}$ from $\tau_-$ to $\tau_+$ produces a path $\gamma$
shown in bold in the Maslov diagram in Figure 2, case (a).  Let
$\Area(\gamma)$ denote the symplectic area enclosed by $\gamma$.  Then
$$ (s_{13}(j_{13}),s_{24}(j_{24}))(\tau_+) = \exp(i \Area(\gamma))
(s_{13}(j_{13},s_{24}(j_{24}))(\tau_-) .$$
\eqref{naive} becomes
$$ | (s_{k,13}(kj_{13}), s_{k,24}(kj_{24}) | \sim \frac{ \cos(\hh \Area(\gamma) +
  \frac{\pi}{4})}{ (k(r-2))^{3/2} \sqrt{\omega( \ppthot,\ppthtf)}} .$$
We claim that 
$$ \Area(\gamma) = 2 \sum_{a < b} l_{ab} \theta_{ab} +
\pi( l_{12} + l_{23} + l_{34} + l_{14}) .$$
As we vary $l_{13}$ and $l_{24}$,
$$ d\Area(\gamma) = 2(\theta_{13} \d l_{13} + \theta_{24} \d l_{24} )
= 2 ( d \sum_{a < b} \theta_{ab} l_{ab} - 2 \Vol(\tau))$$
by Schl\"afli's identity \ref{nonEuclid} \eqref{schlafli}.  Hence
$$ \Area(\gamma) = 2 \sum_{a < b} \theta_{ab} l_{ab} + c(
l_{12},l_{23},l_{34},l_{14}) $$
where $c(l_{12},l_{23},l_{34},l_{14})$ is a constant depending on
$l_{12},l_{23},l_{34},l_{14}$.  For $l_{13},l_{24}$ corresponding to a
degenerate tetrahedron with four exterior edges, we have
$$\theta_{a(a+1)} = \pi,\ \ a = 1,2,3,4, \ \ \ \theta_{13} =
\theta_{24} = 0 .$$
This implies
$$c (l_{12},l_{23},l_{34},l_{14}) = 2\pi( l_{12} + l_{23}
+ l_{34} + l_{14}) .$$ 
Putting everything together, this gives in the transversal
case the naive prediction
$$ |(s_{k,13}(kj_{13}),s_{k,24}(kj_{24}))| \sim \frac{2\pi \cos \left( \frac{k (r-2)}{2\pi} ( \sum
\theta_{ab} l_{ab} - 2  \Vol(\tau)) + \frac{\pi}{4} \right)}
{(k(r-2))^{3/2} \det(\cos(l_{ab}))^{1/4}} .$$
Including half-forms in the quantization scheme should produce shifts
giving Theorem \ref{main} (a).  (Readers familiar with the
representation theory of affine Lie algebras will recognize the
numbers $\hh,2$ in the formulas for $l_{ab}$ and $r(k)$ as the
half-sum of positive weights, resp. dual Coxeter number of
$\lie{sl}_2$.)  Unfortunately we lack the techniques to carry this
out.  In the follow section we introduce a recursion which reduces the
formula to the $q=1$ case; the recursion is needed anyway to handle
the first degenerate case, so even if a geometric proof could be found
for case (a) the proof of the Theorem \ref{main} we give is presumably
shorter.

\section{A non-degenerate tetrahedron}

First we show that both sides satisfy a second order difference
equation as one label is varied; the same strategy was followed by
Mizoguchi-Tada \cite{mi:th} without knowledge of the denominator in
the asymptotic formula.  This implies that each side is a linear
combination of the two linearly independent solutions.  To prove that
the coefficients are equal we show that in the limit as the labels
$j_{ab}$ are made small in relation to the level $r -2$ our equation
reduces to the Ponzano-Regge formula \eqref{pr}.

\subsection{The Schulten-Gordon recursion}

\begin{proposition}\label{recur}
\begin{multline} \label{recureq}
P(j_{23})\left\{
\begin{array}{ccc} j_{12} & j_{23} & j_{13} \\ j_{34} & j_{14} & j_{24} \end{array} \right\}
+[2j_{23}]Q(j_{23}+1)
\left\{
\begin{array}{ccc} j_{12} & j_{23}+1 & j_{13} \\ j_{34} & j_{14} & j_{24} \end{array} \right\}
\\ +[2j_{23}+2]Q(j_{23})
\left\{
\begin{array}{ccc} j_{12} & j_{23}-1 & j_{13} \\ j_{34} & j_{14} & j_{24} \end{array} \right\}= 0,\end{multline}
where 
\begin{align*}
P(&j_{23}) =[2j_{23}][2j_{23}+1][2j_{23}+2][j_{14}+j_{13}+j_{34}+1][j_{13}+j_{34}-j_{14}]\\
&-[2j_{23}][j_{12}+j_{13}-j_{23}][j_{12}+j_{23}-j_{13}+1][j_{24}+j_{34}-j_{23}][j_{24}+j_{23}-j_{34}+1] \\
&-[2j_{23}+2][j_{12}+j_{23}+j_{13}+1][j_{13}+j_{23}-j_{12}][j_{24}+j_{34}+j_{23}+1][j_{34}+j_{23}-j_{24}]
\end{align*}
and
\begin{align*}
Q(j_{23}) &=([j_{12}+j_{23}-j_{13}][j_{12}+j_{13}-j_{23}+1][j_{24}+j_{34}-j_{23}+1][j_{24}+j_{34}-j_{23}]\\
 &[j_{12}+j_{23}+j_{13}+1][j_{13}+j_{23}-j_{12}][j_{24}+j_{34}+j_{23}+1][j_{34}+j_{23}-j_{24}])^{\frac{1}{2}}.
\end{align*}
\end{proposition}
\noindent The proof is a combination of pentagon identities
\eqref{pentagon}, as in \cite{sch:semi}, and is omitted.

In the limit $r \to \infty$ the Schulten-Gordon equation \eqref{receq}
simplifies dramatically.  Let $\Delta_{ab}$ denote the discrete
Laplacian
$$ (\Delta_{ab} f) (j_{ab}) = f(j_{ab} + 1) - 2 f(j_{ab}) + f(j_{ab} -
1) .$$
Let $\tau$ denote the tetrahedron with lengths $ l_{ab} = (j_{ab} +
\hh)/r$, if it exists, and $\theta_{ab}$ the exterior dihedral angle
around the edge $e_{ab}$.  Set
\[
f(j_{ab})=\left( \frac{|G_+(l)|^{\hh}}{\sin(\theta_{ab})} \right)^{\frac{1}{2}} \left\{ \begin{array}{ccc}
j_{12} & j_{23} & j_{13} \\ j_{34} & j_{14} & j_{24} \end{array}
\right\} . \] 

Let $T_S \subset \R^6$ denote the subset of edge
lengths of possibly degenerate tetrahedra, with all 
faces non-degenerate.  $T_S$ is defined by 
the inequalities
\begin{equation} \label{faces}
 |G_+(\tau(d)) > 0 | , 1 \leq a < b < c \leq 4 \end{equation}
$$ | G_+(\tau)| \ge 0 .$$

\begin{theorem}\label{recursion} Let $K$ be a compact
subset of $T_S$ and $1 \leq a < b \leq 4$.  There exists a constant
$C$ and a second order difference operator $\eps(r)$ with coefficients
bounded by $C r^{-2}$ for $r$ sufficiently large, such that
\begin{equation} \label{receq}
 (\Delta_{ab} + 2 -2\cos \theta_{ab} + \eps(r)) f(j_{ab}) = 0
\end{equation}
on the relative interior of $K$.  \end{theorem}

\begin{remark}  One may extend this recursion to the
subset of $\R^6$ of $6$-tuples $l_{ab}$ satisfying only \eqref{faces}
by defining $\theta_{ab}$ to be the imaginary number given by
\eqref{dgram}.
\end{remark}

\begin{proof} Let $cd$ be the opposite edge of $ab$.  
A somewhat long computation involving standard trigonometric
identities shows that
\begin{equation} \label{asymP}  \left| P(j_{ab}) 
\left(\sin{\frac{\pi}{r}} \right)^{5} -  \frac{1}{4}
\frac{\d}{\d \cos(l_{cd})} |G_+(l)| \right|
< C r^{-2}. \end{equation}
Substituting this into the recursion \eqref{recureq} and using 
the factorization \ref{nonEuclid} \eqref{Gram3} gives
\begin{multline} \label{protoprop}
\frac{1}{4}\sin (l_{ab}) \frac{\d|G_+(l)|}{\d \cos(\lE)} \left\{
\begin{array}{ccc} j_{ac} & j_{ab} & j_{bc} \\ j_{bd} & j_{cd} &
j_{ad} \end{array} \right\} +\frac{1}{4}\sin (l_{ab} -\frac{\pi}{r})\cdot
\\ [|G_+(l_{ac},l_{ab}+\frac{2\pi}{r},l_{bc})||G_+(l_{bd},l_{ab}+ \frac{2\pi}{r},l_{ad})|
|G_+(l_{ac},l_{ab},l_{bc})||G_+(l_{bd},l_{ab},l_{ad})|]^{\frac{1}{4}} \\
\left\{\begin{array}{ccc} j_{ac} & j_{ab}+1 & j_{bc} \\ j_{bd} &
j_{cd} & j_{ad} \end{array} \right\} +\frac{1}{4}\sin (l_{ab}
+\frac{\pi}{r})\cdot \left\{\begin{array}{ccc} j_{ac} & j_{ab}-1 &
j_{bc} \\ j_{bd} & j_{cd} & j_{ad} \end{array} \right\} \\
[|G_+(l_{ac},l_{ab}-\frac{2\pi}{r},l_{bc})||G_+(l_{bd},l_{ab}-\frac{2\pi}{r}, l_{ad})|
|G_+(l_{ac},l_{ab},l_{bc})||G_+(l_{bd},l_{ab},l_{ad})|]^{\frac{1}{4}} = O(r^{-2})
\end{multline}
on the relative interior of $K$.  (Here $G_+(l_{ac},l_{ab},l_{bc})$
etc.  refer to Gram matrices of the faces.)  We divide
\eqref{protoprop} by
$$\left(\sin (l_{ab} - \frac{2\pi}{r})\sin (l_{ab}) \sin (l_{ab} + 
\frac{2\pi}{r})
\right)^{\frac{1}{2}}
[|G_+(l_{ac},l_{ab},l_{bc})G_+(l_{ab},l_{bd},l_{ad})|]^{\frac{1}{4}}.$$
Using the equality 
$${\sin (l_{ab})}{(\sin (l_{ab} -  2\pi/r} )
\sin (l_{ab}) \sin ( l_{ab} +  2\pi/r )
)^{\frac{-1}{2}} = \sin(l_{ab})^{-1/2} (1 + O(r^{-2}))$$ 
and similar identities one obtains the vanishing of
\begin{multline}
-\sin (l_{ab})^{-\frac{1}{2}}
[|G_+(l_{ac},l_{ab},l_{bc})||G_+(l_{ab},l_{bd},l_{ad})|]^{-\frac{1}{4}}
\frac{\d |G_+(l)|}{\d \cos(\lE)}
\left\{ \begin{array}{ccc} j_{ac} & j_{ab} & j_{bc} \\ 
                           j_{bd} & j_{cd} & j_{ad} \end{array} \right\}\\ 
+ \sin (l_{ab} +\frac{2\pi}{r})^{-\frac{1}{2}}
[|G_+(l_{ac},l_{ab}+\frac{2\pi}{r},l_{bc})||G_+(l_{bd},l_{ab}+\frac{2\pi}{r},
l_{ad})|]
^{\frac{1}{4}}
\left\{\begin{array}{ccc} j_{ac} & j_{ab}+1 & j_{bc} \\ 
                          j_{bd} & j_{cd} & j_{ad} \end{array} \right\}\\
+\sin(l_{ab} -\frac{2\pi}{r})^{-\frac{1}{2}}[
G_+(l_{ac},l_{ab}-\frac{2\pi}{r},l_{bc})G_+(l_{bd},l_{ab}-\frac{2\pi}{r},l_{ad})
]^{\frac{1}{4}}
\left\{\begin{array}{ccc} j_{ac} & j_{ab}-1 & j_{bc} \\ 
                          j_{bd} & j_{cd} & j_{ad} \end{array} \right\}  
\end{multline}
up to terms of order $r^{-2}$.  Substituting \eqref{dgram} and
\eqref{vol} into the above gives the result.
\end{proof}

\begin{remark} \label{lb} 
An examination of the proof shows that the constant $C$ depends
linearly on a lower bound for
$(|G_+(l_{ac},l_{ab},l_{bc})||G_+(l_{ab},l_{bd},l_{ad})|)
^{\frac{1}{4}}.$
\end{remark}

\subsection{Second-order recursion for the asymptotic formula}

We now show that the right-hand side of \eqref{conj} is also an
approximate solution to \eqref{receq}.  Set
\begin{equation}\label{inf}
\left\{\begin{array}{ccc} j_{12} & j_{23} & j_{13} \\ 
                          j_{34} & j_{14} & j_{24} \end{array} 
\right\}^{\infty}_{q = \exp(\pi i /r)} := 
\frac{2\pi \cos ( \phi + \frac{\pi}{4})}{r^{\frac{3}{2}}
|G_+(l)|^{\frac{1}{4}}}.
\end{equation}
Let $\Delta_{ab}^r$ be the second-order discrete Laplacian at step
size $2\pi/r$,
$$ (\Delta_{ab}^r f)(l) = f(l - 2\pi/r) - 2 f(l) + f(l +
2\pi/r) .$$
Let $T_S^+$ denote the subset of $\R^6$ defined by the inequalities
\eqref{faces} and $ | G_+(\tau)| \ge 0,$ that is, the set
of edge lengths of non-degenerate tetrahedra. 

\begin{theorem}\label{inf-rec}
Let $K$ be a compact subset of $T_S^+$.  There exists a constant $C$
and a second-order difference operator $\eps(r)$ with coefficients
bounded by $C r^{-2}$ for $r$ sufficiently large such that
$$ ( \Delta_{ab}^r + 2 - 2\cos(\theta_{ab}) + \eps(r) )
\left( \frac{|G_+(l)|^{\frac{1}{2}}}{\sin \theta_{ab}}
\right)^{\frac{1}{2}} 
\left\{\begin{array}{ccc} j_{12} & j_{23} & j_{13} \\ 
                          j_{34} & j_{14} & j_{24} \end{array} \right\}^{\infty}_{q
= \exp(\pi i/r)}  
$$
vanishes on the relative interior of $K$.
\end{theorem}

\begin{proof}  Let $\theta = \theta_{ab}$.   By Taylor's theorem
\begin{equation}\label{taylor1}
\phi(l_{ab}\pm \frac{2\pi}{r})
= \phi(l_{ab})\pm \frac{\partial \phi}{\partial 
l_{ab}}
\frac{2\pi}{r}(l_{ab})
+\frac{1}{2}
\frac{\partial^2 \phi}{\partial {l_{ab}}^2}
\left(\frac{2\pi}{r}\right)^2(l_{ab} + \eps)
\end{equation}
for some $\eps \in [-\frac{2\pi}{r},\frac{2\pi}{r}]$.  Because $K$ is
compact, $\frac{\partial^2 \phi}{\partial {l_{ab}}^2}$ is bounded on an
open neighborhood of $K$.  By \ref{nonEuclid} \eqref{schlafli},
\[\frac{\partial \phi}{\partial l_{ab}}=\frac{r}{2\pi}
\theta, \ \ 
\frac{\partial^2 \phi}{\partial
  {l_{ab}}^2}=\frac{r}{2\pi}\frac{\partial \theta}{\partial l_{ab}}.\]
Thus
\begin{equation}
\cos ( \phi(l_{ab} \pm {2\pi}/{r} )) =
\cos(\phi(l_{ab}) \pm \theta)
 -\frac{\pi}{r}\frac{\partial 
\theta}{\partial
  {l_{ab}}}\sin (\phi(l_{ab}) \pm \theta) +O(r^{-2}).
\end{equation}
Similarly 
\[
\sin^{-1/2} ( \theta(l_{ab} \pm {2\pi}/{r}))= \sin^{-1/2} (\theta (l_{ab})) \mp
\frac{\pi}{r}\sin^{-3/2}(\theta(l_{ab})) \cos (\theta) \frac {\partial
\theta}{\partial l_{ab}} +O(r^{-2}).
\]
Expanding $\cos(\phi(l_{ab}) \pm \theta), \sin(\phi(l_{ab}) \pm
\theta)$ one finds that the $O(1)$ and $O(r^{-1})$ terms cancel, which
completes the proof.
\end{proof}

\begin{remark} \label{ub} An examination of the proof shows that the constant
$C$ depends linearly on an upper bound for $ \frac{\partial^2
\theta_{ab}}{\partial {l_{ab}}^2} .$
\end{remark}

\subsection{Euclidean limit of the $6j$ symbol}

Next we investigate the limit of both sides of \eqref{conj} as the
level is taken to infinity much faster than the labels.

\begin{proposition}  \label{euc1} 
Let $j_{ab} \in [0,(r-2)/2] $ and $\delta_{ab} $ be half-integers, for
$1 \leq a < b \leq 4$, and $p \ge 3$ an integer.  Let $j_{ab}(k)
= kj_{ab} + \delta_{ab}$.  Then 
$$ 
\left\{\begin{array}{ccc} j_{12}(k) & j_{23}(k) & j_{13}(k) \\ 
                          j_{34}(k) & j_{14}(k) & j_{24}(k) \end{array} 
\right\}_{q = \exp(\pi i /r(k^p))} \sim  
\left\{\begin{array}{ccc} j_{12}(k) & j_{23}(k) & j_{13}(k) \\ 
                          j_{34}(k) & j_{14}(k) & j_{24}(k) \end{array} 
\right\}_{q = 1}  .$$
\end{proposition}

The proof relies on the 
\begin{lemma} \label{powerp}  There exist  $C,c > 0$
such that if ${r}/{n} < c$ then $| \ln \left( {[n] !}/{ n!}  \right) | <
  Cn^3 r^{-2} $.
\end{lemma}

\begin{proof}  Applying the logarithm to \eqref{qfac} gives 
$$ \ln[n ]! = \sum_{j=1}^{n} \ln \left(\sin \left( \frac{\pi j}{r} \right)
\right) - \ln \left(\sin \left( \frac{\pi}{r} \right) \right).$$
Note that $ | \sin(x) - x | < \frac{x^3}{6} $ and 
$ \ln (\sin(x)) = \ln(x) + (\sin(x) - x)/y $
for some $y \in [x,\sin(x)]$.  For $x$ sufficiently small, 
$$ \sin(x)
\in [x - \frac{x^3}{6},x], \ \ \ \ y > x/2 $$
which implies
$$ |\ln(\sin(x)) - \ln(x)| < \frac{x^3}{6} \frac{2}{x} = \frac{x^2}{3} .$$
Hence
$$ | \ln [n]! - \ln n! | 
< \sum_{j = 1}^{n} \frac{1}{3} ((\frac{\pi j}{r})^2  + (\frac{\pi}{r})^2) 
<  \frac{n}{3} ((\frac{\pi n}{r})^2 + (\frac{\pi}{r})^2) 
< \frac{2\pi^2 n^3}{3r^2}   .$$
\end{proof}

\begin{proof}[Proof of \ref{euc1}:] 
By Lemma \ref{powerp} and Racah's formula \eqref{racah}, since the
number of terms in the sum grows linearly with $k$.
\end{proof}

\subsection{Euclidean limit of the asymptotic formula}

In this section we compare the quantities 
$$ 
\left\{\begin{array}{ccc} j_{12} & j_{23} & j_{13} \\ 
                          j_{34} & j_{14} & j_{24} \end{array} 
\right\}^{\infty}_{q = \exp(\pi i /r)} := 
\frac{2\pi \cos ( \phi +
  \frac{\pi}{4})}{r^{\frac{3}{2}}|G_+(l)|^\qq}
$$
and
$$ 
\left\{\begin{array}{ccc} j_{12} & j_{23} & j_{13} \\ 
                          j_{34} & j_{14} & j_{24} \end{array} 
\right\}^{\infty}_{q = 1} := 
\frac{2\pi \cos ( \phi_0 + \frac{\pi}{4})}{||G_0(j + \hh)||^{1/4}}
$$
in the case that the corresponding tetrahedra are non-degenerate.  We
will prove

\begin{proposition} \label{euc2}  
Let $j_{ab} \in [0,(r-2)/2], 1 \leq a < b \leq 4 $ be such that a
non-degenerate Euclidean tetrahedron with edge lengths $j_{ab}$
exists.  Let $\delta_{ab} $ be half-integers, and $p > 3/2$.  Let
$j_{ab}(k) = kj_{ab} + \delta_{ab}$.
$$ 
\left\{\begin{array}{ccc} j_{12}(k) & j_{23}(k) & j_{13}(k) \\ 
                          j_{34}(k) & j_{14}(k) & j_{24}(k) \end{array} 
\right\}_{q = \exp(\pi i /r(k^p))}^\infty \sim  
\left\{\begin{array}{ccc} j_{12}(k) & j_{23}(k) & j_{13}(k) \\ 
                          j_{34}(k) & j_{14}(k) & j_{24}(k) \end{array} 
\right\}_{q = 1}^\infty .$$

\end{proposition}
\begin{proof}  By \ref{nonEuclid} 
\eqref{GG0} $r(k^p)^{3/2} |G_+(l)| \sim ||G_0(j + \hh)||$.  Since the
edge lengths are $O(k^{1-p})$, the volume $\Vol(\tau)$ is
$O(k^{3-3p})$.  Let $\theta_{ab}^E$ denote the dihedral angles of the
Euclidean tetrahedron with edge lengths $l_{ab}$, if it exists.  It
follows from Corollary \ref{angle} that $ \cos(\theta_{ab}) =
\cos(\theta_{ab}^E)(1 + O(k^{2 - 2p})) .$ Since $\theta_{ab} \to
\theta_{ab}^E$ as $k \to \infty$, $ \theta_{ab} - \theta_{ab}^E < C
k^{2 - 2p} .$ Hence
$$ \sum l_{ab} (\theta_{ab} - \theta_{ab}^E)
< Ck^{3 - 3p} $$
which implies that 
$$ \phi(k) - \phi_0(k) < C k^{3 - 2p} .$$
For $p > 3/2$, this approaches zero as $k \to \infty$. 
\end{proof}

\subsection{All edge lengths equal}

The second order difference equation corresponds to a system of
first-order difference equations in the standard way.  In fact, if $v$
is the vector formed from all $2^6$ sets of labels of the form $j_{ab}
+ \delta_{ab}$, where $\delta_{ab} = 0$ or $1$, then the vector of the
corresponding $2^6$ $6j$ symbols satisfies a first order recursion in
all six variables $j_{ab}$.  We call the first-order systems
corresponding to \eqref{recureq},\eqref{receq} the first-order
Schulten-Gordon difference equation and first-order asymptotic
difference equation, respectively.  By Theorems \ref{recursion} and
\ref{inf-rec}, to show Theorem \ref{main} (a) it suffices to prove the
claim for all components of $v$ and a single set of values of
$j_{ab}$.  We do so for the case that all edge lengths are equal.  By
\ref{nonEuclid} \eqref{equal}, a spherical tetrahedron with edge
lengths equal to $l$ exists for $0 \leq l \leq \arccos( -1/3)$.  Let
$h = h(k)$ be a function of $k$ so that $h > k^p$, for some $p > 3$.
Choose $r >> 0$ so that there exists $j \in
(0,(r-2)/4) \cap \Z/2$.  For sufficiently large $k$, there exists an
approximation $\gamma_k(t)$ to the linear path $lt, t \in [ k/h(k),1]$
such that
\begin{enumerate}
\item $\gamma_k(t) + \delta_{ab}/r(h(k))$ is the set of edge lengths
of a non-degenerate spherical tetrahedron, for any choice of
$\delta_{ab}$;
\item $\gamma_k$ is a lattice path for $\Z/2r(h(k))$, that is,
obtained by concatenating paths between lattice points whose
difference is a standard basis vector divided by $2r(h(k))$;
\item $ \gamma_k(k/h(k))_{ab} = (kj +
\hh)/h(k)$ and $\gamma_k(1)_{ab} = (h(k)j + \hh)/h(k)$;
\item the number of lattice points in $\gamma_k$ is $O(h(k))$.
\end{enumerate}
Let $l(k,m), m=0 ,\ldots, n(k)$ be the lattice points occurring in
$\gamma$.  Let $A_m$ be the operator in the $m$-th term of the
first-order Schulten-Gordon difference equation, and $A_\m^\infty$ the
operator in the $m$-th term in the asymptotic difference equation.  
On an open cone containing $l$, we have an estimate
$$(|G_+(l_{ac},l_{ab},l_{bc})|
|G_+(l_{ab},l_{bd},l_{ad})|)^{\frac{1}{4}} > C \Vert l \Vert^{2}.$$
By Remark \ref{lb} $A_m^{-1} A_m^\infty = I + \eps(m)$ where $\eps(m)$
has coefficients bounded by
\begin{equation} \label{error} C h^{-2} \Vert l(k,m) \Vert^{-2} .\end{equation}
Let 
$$v(k,m) = A(k,m) A(k,m-1) \ldots A(k,1) v(k,0) $$
denote the solution to the Schulten-Gordon difference equation, and 
$$v^\infty(m) = A^\infty(m) A^\infty(m-1) \ldots A^\infty(1) v(0) $$ 
the solution to the asymptotic equation. Then 
\begin{eqnarray*}
 v^\infty(m) &=& A(m) (I + \eps(m)) A(m-1) (I + \eps(m-1)) \ldots
A(1) (I + \eps(1)) v(0) \\
&=& (I + \eps'(m))(I + \eps'(m-1)) \ldots (I + \eps(1))
v(m) \end{eqnarray*}
where $\eps'(1),\ldots,\eps'(m)$ are operators with coefficients also
bounded by \eqref{error}.  The norm of the product of error factors
satisfies the estimate
$$ \left\Vert \ln\left( \prod_{m=1}^{n(k)} (I + \eps'(m)) \right) 
\right\Vert \leq C \sum_{m= 1}^{n(k)} \Vert h l(k,m)
\Vert^{-2}  $$
which goes to zero as $ k \to \infty$.  It follows that the components
of $v^\infty(m)$ are asymptotic to those of $v(m)$ as $k \to \infty$.
(The convergence of approximate solutions of second order difference
equations to actual solutions is discussed in much greater generality
in \cite{ge:tp}.)

A similar estimate holds for the right-hand side of \eqref{conj}.  On
the cone generated by a small neighborhood of $l_{ab}^0$ we have a
bound
$$ \frac{\partial^2 \theta_{l_{ab}}}{\partial {l_{ab}}^2} (l) < C \Vert l
\Vert^{-2}.$$
Using Remark \ref{ub}, the error factor approaches $1$ as $k \to
\infty$, that is, the asymptotic $6j$ symbols approach an exact
solution of the asymptotic difference equation.

To show that $v_k(n(k))$ is asymptotic to $v_k^\infty(n(k))$ it
remains to show that each component of $v_k(0)$ is asymptotic to the
corresponding component of $v_k^\infty(0)$.  By Propositions
\ref{euc1}, \ref{euc2} it suffices to show that for any $\delta_{ab}
\in \{ 0, 1 \}$,
$$ \bsj kj + \delta_{12} & k j + \delta_{23} & kj + \delta_{13} \\ 
        k j+ \delta_{34}  & kj + \delta_{14}  & kj + \delta_{24} \esj_{ q = 1}
\sim 
 \bsj kj + \delta_{12} & k j + \delta_{23} & kj + \delta_{13} \\ 
        k j+ \delta_{34}  & kj + \delta_{14}  & kj + \delta_{24} \esj_{ q = 1}^\infty .$$
For $\delta_{ab} = 0$ this is Roberts theorem \cite{ro:qu}.
To handle the remaining cases, we need a generalization:

\begin{theorem} \label{gen}
Let $j_{12}(k),\ldots,j_{34}(k) \in \Z/2$ be functions of $k$ such
that $l_{ab}(k) = 2\pi (j_{12}(k) + \hh)/k$ converge as $k \to \infty$
to $l_{ab} \in (0,\pi)$ the edge lengths of a non-degenerate Euclidean
tetrahedron.
$$ 
\left\{\begin{array}{ccc} j_{12}(k) & j_{23}(k) & j_{13}(k) \\ 
                          j_{34}(k) & j_{14}(k) & j_{24}(k) \end{array} 
\right\}_{q = 1} \sim  
\left\{\begin{array}{ccc} j_{12}(k) & j_{23}(k) & j_{13}(k) \\ 
                          j_{34}(k) & j_{14}(k) & j_{24}(k) \end{array} 
\right\}_{q = 1}^\infty  .$$
\end{theorem}

\begin{proof}   It seems that this follows by the same method used by 
Roberts \cite{ro:qu}.  Alternatively, the $6j$ symbols for $q = 1$ as
a function of any $j_{ab}$ satisfy the Schulten-Gordon second-order
difference equation up to an error term depending linearly on an upper
bound for $ ||G_0(l_{ab},l_{ac},l_{bc})||
G_0(l_{ab},l_{ad},l_{bd})||^{-1/4}$.  Let $T_E^+$ be the subset of
$\R^6$ consisting of edge lengths of non-degenerate tetrahedra, and $K
\subset T_E^+$ a compact subset.  On $k K$, the $6j$ symbols solve the
difference equation up to an error operator of order $O(k^{-2})$, and
so equal an exact solution up to $O(k^{-1})$.  Any exact solution is a
linear combination
$$2\pi r^{-3/2} ||G_0(l)||^{-1/4} (c_1 \cos(\phi_0 + \pi/4) + c_2
\sin(\phi_0 + \pi/4)) ;$$
We wish to show that $c_1 = 1, c_2 = 0$.  Let $ab = 12$ and suppose
that $j_{ab}(k) = kj_{ab}$ for $ab \neq 12$.  Roberts' theorem applied
to $\left\{\begin{array}{ccc} kj_{12} & kj_{23} & kj_{13} \\ kj_{34} &
kj_{14} & kj_{24} \end{array} \right\} $ gives
\begin{equation} \label{one}
 c_1 \cos(\phi_0 + \pi/4) + c_2 \sin(\phi_0 + \pi/4) = \cos(\phi_0 + \pi/4). 
\end{equation}
Roberts' theorem applied to $\left\{\begin{array}{ccc} k(j_{12}+1) &
kj_{23} & kj_{13} \\ kj_{34} & kj_{14} & kj_{24} 
\end{array} \right \}$ gives
\begin{equation}  \label{two}
 c_1 \cos(\phi_0 + \theta_{12} + \pi/4) + c_2 
\sin(\phi_0 + \theta_{12} + \pi/4) = \cos(\phi_0 + \theta_{12} + \pi/4) .\end{equation}
Since $\theta_{12} \in (0,\pi)$, the equations \eqref{one},\eqref{two}
are linearly independent, hence $c_1 = 1 $ and $c_2 = 0$.  Applying
the same argument to the variation in $j_{23}, j_{13}$ etc. proves the
theorem.
\end{proof}
\noindent This completes the proof of Theorem \ref{main} (a).  The
argument used for Theorem \ref{gen} applies to the $6j$ symbols for $
q = \exp(\pi i/ r)$ and gives the following generalization of Theorem
\ref{main} (a).
\begin{theorem}  \label{gen1}
Let $r >2, j_{12}(k),\ldots,j_{34}(k)
\in [0,(r-2)/2] \cap \Z/2$ be functions of $k$ such that 
$l_{ab}(k) = 2\pi (j_{12}(k) + \hh)/r(k)$ converge as $k \to \infty$ to 
$l_{ab} \in (0,\pi)$ the edge lengths of a non-degenerate
spherical tetrahedron. 
$$ 
\bsj j_{12}(k) & j_{23}(k) &
 j_{13}(k) \\ j_{34}(k) & j_{14}(k) & j_{24}(k) \esj_{q = \exp(\pi i /r(k))} 
 \sim 
\bsj j_{12}(k) & j_{23}(k) &
 j_{13}(k) \\ j_{34}(k) & j_{14}(k) & j_{24}(k) \esj_{q = \exp(\pi i/r(k))}^\infty .$$
\end{theorem}

\section{A degenerate tetrahedron with non-degenerate faces}

We apply the difference equation of the previous section to a sequence
of labels in which the corresponding tetrahedra degenerates to a
degenerate tetrahedron with non-degenerate faces.  The two possible
limiting cases are shown below in Figure 6.  Let
$l_{ab}^0,\theta_{ab}^0$ denote the lengths, resp.  angles of the
degeneration.
\begin{figure}[h]
\setlength{\unitlength}{0.00033333in}
\begingroup\makeatletter\ifx\SetFigFont\undefined%
\gdef\SetFigFont#1#2#3#4#5{%
  \reset@font\fontsize{#1}{#2pt}%
  \fontfamily{#3}\fontseries{#4}\fontshape{#5}%
  \selectfont}%
\fi\endgroup%
{\renewcommand{\dashlinestretch}{30}
\begin{picture}(7524,2739)(0,-10)
\path(912,2712)(12,912)(2112,12)
	(2412,1812)(912,2712)(2112,12)
\path(2412,1812)(12,912)
\path(5712,2712)(4512,312)(7512,612)
\path(5712,2712)(7512,612)(6012,1212)
\path(5712,2712)(6012,1212)
\path(6012,1212)(4512,312)
\end{picture}
}
\caption{Two limiting cases}
\end{figure}
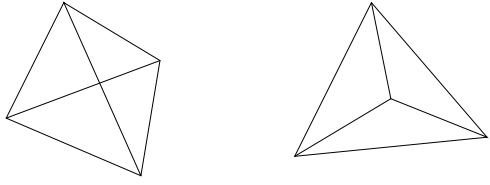
\noindent Following Schulten and Gordon \cite{sch:semi} we define
$$ \phi_0 = \frac{r(k)}{2\pi} \sum_{a <b} \theta_{ab}^0 l_{ab} 
 \ \ \ \ \  \theta_{ab}^0 = 
\begin{cases}
         0,  &\text{if $\theta_{ab} \leq \pi/2$\ ;} \\
         \pi,  &\text{if $\theta_{ab} > \pi/2$\ .} 
\end{cases}
$$
In the two limiting cases shown above, we have 
$$ \phi_0 = \pi ( j_{12} + j_{23} + j_{34} + j_{41} + 2) , 
\  \ \text{resp.} \  \pi (j_{12} + j_{23} + j_{31} + \thh). $$
Hence $ \phi_0$ is an integer, resp. half integer, times $\pi$.  By
Sch\"afli's formula \ref{nonEuclid} \eqref{schlafli}
\begin{equation} \label{phio}
\d (\phi - \phi_0) = \frac{r(k)}{2\pi} \sum (\theta_{ab} - \theta_{ab}^0) \d \l_{ab}
\end{equation}
which implies $ \phi - \phi_0 < 0 \text{,} \ \text{resp.} \ \phi -
\phi_0 > 0 $ for tetrahedra near the first, resp. second, limiting
case.

We apply the LG-WKB method to arrive at a conjectural formula.  In
order for $f$ to be a solution to the recursion relation
\eqref{receq}, $f_1,f_2$ must satisfy
$$ (\Delta_{ab} + 2 - 2 \cos(\theta_{ab} - \theta_{ab}^0) )
f(j_{ab}) = 0 .$$
According to ansatz described by Schulten and Gordon, in the
semiclassical limit the solution must solve the differential equation
\begin{equation} \label{ansatz}
 \left( \frac{\partial^2}{\partial j_{ab}} - (\theta_{ab} -
\theta_{ab}^0)^2 \right)
\sqrt{ \frac{ \sin(\theta_{ab} - \theta_{ab}^0) }{\theta_{ab} -
    \theta_{ab}^0}}
f(j_{ab}) = 0 .\end{equation}
Let $F(x)$ denote a solution to the Airy equation $ F''(x) = x F(x) .$
Define
$$ f(j_{ab}) = \left(\frac{ \sin(\theta_{ab} - \theta_{ab}^0) }{\theta_{ab} -
    \theta_{ab}^0} \right)^{-1/2} A(j_{ab}) F(\Omega(j_{ab})) .$$
Inserting this into \eqref{ansatz} gives
$$ \left( \frac{A''}{A} + \Omega'^2 \Omega + (\theta_{ab} - \theta_{ab}^0)^2 \right)
F(\Omega(j_{ab})) + \left( 2 \frac{A'}{A} \Omega' + \Omega''\right)
F'(\Omega(j_{ab})) = 0 .$$
Assuming that $A''/A \approx 0$, the solution is
$$ \Omega(j_{ab}) = -\frac{3}{2} \left(\int (\theta_{ab} - 
\theta_{ab}^0) \d j_{ab} \right)^{2/3} = -
\frac{3}{2} (\phi(j_{ab}) - \phi^0(j_{ab}))^{2/3} $$
$$ A(j_{ab}) = C | \phi(j_{ab})
- \phi^0(j_{ab}) |^{1/4} \theta_{ab}^{-1/2} .$$
We need to choose $F,f$ to agree with the solution in the transversal
case.  Introduce the regular and irregular Airy functions
$$ \Ai(x) = \frac{1}{\pi} \int_0^\infty
\cos(t^3/3 + xt) \d t , 
\ \ \ \Bi(x) = \frac{1}{\pi} \int_0^\infty
(\exp(-t^3/3 + xt) + \sin(t^3/3 + xt) ) \d t .$$ 
For large $x$ one has the asymptotics 
\begin{equation} \label{Airyapprox}
 \Ai(-x) \sim \pi^{-1/2} x^{-1/4} \cos(\xi - \frac{\pi}{4}), 
\ \ \ 
\Bi(-x) \sim \pi^{-1/2} x^{-1/4} ( - \sin(\xi -
\frac{\pi}{4})) 
\end{equation}
where $\xi = \frac{2}{3} x^{3/2} .$     It follows that the function
\begin{equation} \label{uniform}
 f =  \frac{2 \pi^{3/2} Z^{1/4}}{r^{3/2}\sin(\theta_{ab})^{1/2}}
 \begin{cases}
\cos(\phi_0) \Ai(-Z) - \sin(\phi_0) \Bi(-Z) & \phi - \phi_0 < 0 \\
\cos(\phi_0) \Bi(-Z) - \sin(\phi_0) \Ai(-Z) & \phi - \phi_0 < 0 
\end{cases}
\end{equation}
where 
$$ Z =  (\frac{3}{2} | \phi - \phi_0 |) ^{2/3} $$
satisfies the asymptotic recursion and matches up with the
non-degenerate solution.  

\begin{lemma} \label{limlem}  
%
In the limit $\theta_{ab} \to \theta_{ab}^0$,
$$ \frac{\phi - \phi_0}{|G_+(l)|^{3/2}} \to \frac{3r}{2\pi A_1A_2A_3A_4}.$$
%
%
%
\end{lemma}

\begin{proof}  
Using \eqref{vol}, at $l_{ab} = l_{ab}^0$ we have
\begin{eqnarray*}
\pplab \frac{|G_+(l)|^{3/2}}{A_1A_2A_3A_4}
&=&  \frac{3}{A_1A_2A_3A_4} |G_+(l)| 
 \pplab |G_+(l)|^{1/2} \\
&=&  \frac{3}{A_1A_2A_3A_4} 
|G_+(l)|^{1/2} 
\frac{A_a A_b  \sin(\theta_{ab})}{ \sin(l_{ab})}
\pplab \frac{ A_c A_d  \sin(\theta_{cd})}{ \sin(l_{cd})} \\
&=&  \frac{3 |G_+(l)|^{1/2} }{A_1A_2A_3A_4} 
\frac{A_a A_b  \sin(\theta_{ab})}{\sin(l_{ab})}
\frac{A_c A_d \cos(\theta_{cd})}{\sin(l_{cd})} 
 \frac{\partial  \theta_{cd}}{\partial l_{ab}} \\
&=& 3 |G_+(l)|^{1/2} 
\frac{\sin(\theta_{ab}) \cos(\theta_{cd})}{\sin(l_{ab}) \sin(l_{cd})}
\frac{\sin(l_{ab})\sin(l_{cd})}{|G_+(l)|^{1/2}} \\
&= & 3   \cos(\theta_{cd}) \sin(\theta_{ab}) .
\end{eqnarray*}
The Lemma now follows from L'Hopital's rule.
%
%
\end{proof}

 The values of the Airy functions at zero are
$$ \Ai(0) = (3^{2/3} \Gamma(\frac{2}{3}))^{-1}, \ \ \Bi(0) = (3^{\frac{1}{6}} \Gamma(\frac{2}{3}))^{-1}
.$$
Therefore, the limit of the $6j$ symbol as $j_{ab} \to j_{ab}^0$ is
$(-1)^{j_{12} + j_{23} + j_{34} + j_{41}} $ resp.  $(-1)^{j_{12} +
j_{23} + j_{13}}$ times
\begin{eqnarray*}
\lim_{j_{ab} \to j_{ab}^0} \frac{\sin(\theta_{ab})^{1/2}}{| G_+(l)|^{1/4}}
f(j_{ab}) &=& \lim_{j_{ab} \to j_{ab}^0} \frac{2 \pi^{3/2} 
((\frac{3}{2} |\phi - \phi_0|)^{2/3})^{1/4}}{r^{3/2} | G_+(l)|^{1/4}}
\Ai(0)  \\
&=&  2\pi^{3/2} (3/2)^{1/6} \left(  \frac{3r}{2\pi A_a A_b A_c
  A_d}   \right)^{\frac{1}{6}}    3^{-2/3} \Gamma \left(\frac{2}{3}\right)^{-1} \\
&=&  r^{-4/3} 2^{2/3} 3^{-1/3} \pi^{4/3} (A_a A_b A_c A_d)^{-1/6} \Gamma \left(\frac{2}{3}\right)^{-1}
\end{eqnarray*}
which is (b) in Theorem \ref{main}.

The LG-WKB method for finite difference equations near a turning point
is studied by Geronimo, Bruno and Assche in \cite[Theorem 3.8]{ge:tp}.
The assumptions \cite[p.106]{ge:tp} that guarantee that an approximate
solution converges to the solution given by the ansatz translate to
the conditions
\begin{equation} \label{i}
\left( \frac{\partial}{\partial l_{ab}} \right)^i (\theta_{ab} -
\theta_{ab}^0)^2 \in C^0((0,\pi)), \ i=0,1,2,3; \end{equation}
\begin{equation} \label{ii}
\left( \frac{\partial}{\partial l_{ab}} \right)^i \frac{(\theta_{ab} -
\theta_{ab}^0)^2}{l_{ab} - l_{ab}^0} \in C^0((0,\pi)), \ i = 0,1,2.
\end{equation}
To see that these conditions holds, note that by \eqref{dgram}
$\cos^2(\theta_{ab} - \theta_{ab}^0)$ is a smooth function of
$l_{ab}$.  Hence so is $(\theta_{ab} - \theta_{ab}^0)^2$ for
$\theta_{ab}$ sufficiently close to $0,\pi$, which proves \eqref{i}.
\eqref{ii} follows since $(\theta_{ab} - \theta_{ab}^0)^2$ has a
simple zero at $l_{ab} = l_{ab}^0$.  Condition (iii) on
\cite[p.106]{ge:tp} is achieved by choosing the sign of the parameter
$t$ appropriately, for $t = \pm(l_{ab} - l_{ab}^0)$ sufficiently
small.

Unfortunately we do not know whether there are any configurations of
the type described in Theorem \ref{main} (b), that is, degenerate
tetrahedra with non-degenerate faces and all edge lengths commensurate
with $\pi$.  However, the following generalization of Theorem
\ref{main} (b) is non-vacuous:
\begin{theorem}  \label{gen2}  Let $j_{ab}(k) \in [0,(r-2)/2] \cap \Z/2$
be a sequence of labels such that $l_{ab}(k) = (j_{ab}(k) + \hh)/r(k)$
converge to the edge lengths $l_{ab}$ of a degenerate tetrahedron with
non-degenerate faces and $l_{ab}(k) - l_{ab} = O(k^{-1})$.  The
sequence $\bsj j_{12}(k) & j_{23}(k) & j_{13}(k) \\ j_{34}(k) &
j_{14}(k) & j_{24}(k) \esj $ is asymptotic to the expression in
Theorem \ref{main} (b) as $ k\to \infty$.
\end{theorem}
\begin{proof}  By the discussion above, the sequence has asymptotics
  given by \eqref{uniform}.  By Schl\"afli's formula and Taylor's
theorem, $\phi - \phi^0 < C(\theta_{ab} - \theta_{ab}^0) = O(k^{-1/2})
.$ Hence the argument in the Airy function in \eqref{uniform} goes to
zero, which completes the proof.
\end{proof} 

\section{One face degenerate} 

If one face, say $123$, is degenerate then there are two
possibilities: either one lengths, say $j_{13}$, is the sum of the
other two, $j_{12} + j_{23}$, or $j_{12} + j_{23} + j_{23} = r-2$.  In
the first case the Racah sum \eqref{racah} has a single term
\begin{multline} 
 \bsj j_{12} & j_{23} & j_{12} + j_{23} \\
        j_{34} & j_{14} & j_{24} \esj
=  
\left\{
\frac{
[2j_{12}]!
[2j_{23}]!
[j_{12} + j_{23} + j_{34} + j_{14} + 1]!
}{
[2j_{12} + 2j_{23} + 1]!
[-j_{12} - j_{23} + j_{34} + j_{14}]!
} \right\}^{\hh} \\
(-1)^{j_{12} + j_{23} + j_{34} + j_{14}} \left\{
\frac{
[j_{12} + j_{23} + j_{34} - j_{14} ]!
[j_{12}+j_{23} - j_{34} + j_{14}]!}
{[j_{23} + j_{34} - j_{24}]!
[j_{23} - j_{34} + j_{24}]! 
[j_{23} + j_{34} + j_{24} + 1]!
} \right\}^{\hh} \\ 
\left\{ \frac{
[j_{34} + j_{24} - j_{23}]!
[j_{14} + j_{24} - j_{12}]!}
{
[j_{12} + j_{24} + j_{14} + 1]! 
[j_{12} - j_{14} + j_{24}]!
[j_{12} + j_{14} - j_{24}]! 
}
\right\}^{\hh} .
\end{multline}
The asymptotics of the $6j$ symbol are therefore determined by the
asymptotics of quantum factorials which have been investigated by Moak
\cite{mo:qa}, see also \cite{ue:qg}.  The following can be derived
from \cite[(2.12-16)]{mo:qa}.

\begin{proposition}  (q-Stirling) 
Set $q = \exp(\pi i /kr)$.  Then as $k \to \infty$,
\begin{equation} \label{qStirling}
 {[kn]!} \sim \sqrt{2 \pi} [kn]^{kn+\hh} I(n\pi/r)^{-kn} \end{equation}
where 
$$ I(x) := \sin(x) \exp \left(
- x^{-1} \int_0^{x} \ln(\sin(y)) \d y \right),
\ \ \lim_{x \to 0} I(x) = e.$$
\end{proposition}

We conjecture that the factors involving $I(n\pi/r)$ cancel; in the
case $q = 1$ these factors are replaced by factors of $e$ which do
cancel.  Assuming this, one obtains
\begin{equation} \label{sk}
 s(k) \sim 2^{-1/4} r(k)^{-\frac{5}{4}}\pi \left\{
\frac{\sin(l_{12})^{\hh} \sin(l_{23})^{\hh}I(l_{13},l_{34},l_{41})}
{\sin(l_{12} + l_{23})^{\thh}
I(l_{12},l_{24},l_{41})I(l_{23},l_{34},l_{41})} \right\}^{1/2}
\end{equation}
times $ (-1)^{k(j_{13} + j_{34} + j_{14})}$ where
$$I(l_{ab},l_{bc},l_{ac}) = 
\frac{H(l_{ab},l_{bc},l_{ac})^{1/2}
  \sin(\hh (l_{ab}+l_{bc}+l_{ac}))}{\sin(\hh (l_{bc}+l_{ac}-l_{ab}))} $$
and
$$ H(l_{ab},l_{bc},l_{ac}) =
\begin{array}{c} \sin(\frac{l_{ab}+l_{bc}+l_{ac}}{2}) \sin(\frac{l_{ab}+l_{bc}-l_{ac}}{2})
\sin(\frac{l_{ab}-l_{bc}+l_{ac}}{2})\sin(\frac{-l_{ab}+l_{bc}+l_{ac}}{2})
\end{array} .$$
If $j_{12} = j_{23}$, \eqref{sk} simplifies to
$$ s(k) \sim \frac{r(k)^{-5/4} \pi}{\sin(l_{12})^{3/4}
 \sin(l_{24})^{1/2}} .$$
The analogous formula for $q = 1$ was checked numerically.

We do not know whether there exist configurations of the type
described in Figure 2(c) in the spherical case with edge lengths
commensurable to and less than $\pi$.  Neither to do we have a uniform
formula, similar to \eqref{uniform}, which describes the asymptotics
of the $6j$ symbols in the neighborhood of such a degeneration.  The
case $j_{12} + j_{23} + j_{13} = r-2$ seems even more
mysterious.

\section{Colinear vertices or degenerate edges} 

Suppose that every face of a tetrahedron $\tau$ is degenerate.  We may
assume without loss of generality that the vertices lie on a geodesic
with length at most $\pi$ in the order $1,2,3,4$ so that $ j_{ab} +
j_{bc} = j_{ac} $ for any $1 \leq a < b < c \leq 4$.  In this case the
$6j$ symbol equals
$$ (-1)^{2kj_{14}} ([2j_{13} + 1] [2j_{24} + 1])^{-1/2}
\sim (-1)^{2kj_{14}} \pi r(k)^{-1} (\sin(l_{13}) \sin(l_{24}))^{-1/2} .$$

In the case that one of the lengths is zero one obtains from
\eqref{racah} the formula
\begin{eqnarray*}
 \bsj kj_{12} & kj_{23} & kj_{13} \\
     kj_{34} & kj_{12} & 0 \esj &=&
(-1)^{kj_{12}+ kj_{23} + kj_{13}} ([2kj_{12}+1] [2kj_{23}+1])^{-1/2} \\
&\sim& (-1)^{k(j_{12} + j_{23} + j_{13})} \pi r(k)^{-1}
\sin(l_{12})^{-\hh} \sin(l_{23})^{-\hh} \end{eqnarray*}
if both lengths are non-zero.  Similarly, one obtains the last two
formulas in Theorem \ref{main}.

\section{The classically forbidden case}

Naturally one expects $6j$ symbols are exponentially decreasing in $k$
if $\tau$ does not exist; this was proved that in the case $q = 1$.
The difference equation of Section 4 together with the results of
Geronimo, Bruno and Assche \cite{ge:tp} prove that the $6j$ symbols
are exponentially decreasing as a function of a label $j_{ab}$ after
it passes into the classically forbidden regime. However, it is not
clear to us whether this implies exponential decay as a function of
$k$ for arbitrary classically forbidden tetrahedra.

\section{The hyperbolic case}

In the case $q = \exp(\pi / r)$ (no factor of $i$) the asymptotics are
related to hyperbolic tetrahedra.  If it exists, let $\tau(k)$ denote
the hyperbolic tetrahedron with edge lengths $ l_{ab}(k) = 2\pi
(kj_{ab} + \hh)/r(k) $ and exterior dihedral angles $\theta_{ab}(k)$.
Let $\tau$ denote the hyperbolic tetrahedron with edge lengths
$ l_{ab} = 2\pi j_{ab}/(r-2) $
and exterior dihedral angles $\theta_{ab}$.  One has identities
similar to the spherical case, e.g.
$$ |G_-(\tau)|^{\hh} = \frac{A_a A_b \sin(\theta_{cd})}{\sinh(l_{cd})} .$$
where $ A_a = |\det(-\cosh(l_{bc} |_{b,c \neq a}))|^{1/2} .$ We state
the results in the hyperbolic case as follows.  They are identical to
the spherical case, except for the replacements $\sin$ with $\sinh$,
$\cos$ with $-\cosh$, and $\Vol(\tau)$ with $-\Vol(\tau)$).  The
proofs are similar.

\begin{theorem} \label{mainh} Let $j_{ab}, 1 \leq a < b \leq 4$
be non-negative half-integers, $r > 2$ and
$$s(k) := \bsj kj_{12} & kj_{23} & kj_{13} \\ kj_{34} & kj_{14} &
 kj_{24} \esj_{q = \exp(\pi /r(k))} .$$
\begin{enumerate}
\item If $\tau$ exists and is non-degenerate, then
$$ s(k) \sim \frac{2\pi \cos( \phi(k) + \pi/4)}{r(k)^{3/2}
  |\det(-\cosh(l_{ab}))|^{1/4}} $$
where
$$\phi(k) = \frac{r(k)}{2\pi} \left(\sum_{a < b} \theta_{ab}(k) l_{ab}(k) +
2\Vol(\tau(k)) \right).$$
\item If $\tau$ exists, has zero volume but all faces have non-zero
area then 
$$
s(k) \sim r(k)^{-\frac{4}{3}} \pi^{\frac{4}{3}} 2^{\frac{2}{3}} 3^{-\frac{1}{3}}
\Gamma(\frac{2}{3})^{-1} 
\frac{\cos(k \sum \theta_{ab} j_{ab})}{(A_1 A_2
A_3 A_4)^{1/6}} .$$
%
%
%
\item[(d)] If $\tau$ exists and has exactly one edge length, say
$l_{ab}$ vanishing, then 
$$
s(k) \sim  (-1)^{k(j_{bc} + j_{cd} + j_{bd})}
\pi  r(k)^{-1} ( \sinh(l_{ac}) \sinh(l_{bd}))^{-1/2} .$$
\item[(e)] If $\tau$ exists, all faces have zero area but all edge
lengths are non-vanishing then $\tau$ lies on a geodesic, say
with vertices in order $a,b,c,d$ then 
$$ 
s(k) \sim  (-1)^{2kj_{ad}} \pi r(k)^{-1} (\sinh(l_{ac})
\sinh(l_{bd}))^{-1/2}.$$
\item[(f)] If $\tau$ exists, all faces are degenerate, and all non-zero
edge lengths are equal then (supposing without loss of generality that
$l_{ad} \neq 0$) 
$$ s(k) \sim  
(-1)^{2kj_{ad}}
(\frac{\pi}{r(k)})^{1/2} \sinh(l_{ad})^{-1/2} .$$
\item[(g)] If $\tau$ exists, but all edge lengths are vanishing
then $ s(k) = 1 .$
\end{enumerate}
\end{theorem}

\section{Questions}

\begin{enumerate}
\item Is there a uniform formula which includes all cases, just as the
Schulten-Gordon formula \eqref{uniform} includes cases (a) and (b)?
Such a formula (or least uniform estimates) would be needed to apply
these results to the asymptotics of the quantum invariants such as
Turaev-Viro and Jones, since one would presumably have to show that
the contributions from the non-degenerate case are dominant.
\item Is there a geometric description of the tensor category of
representations of $U_q(\lie{sl}_2)$ in terms of the moduli space of
hyperbolic bundles \cite{al:hy}, for $q$ positive real, similar to the
description for $q$ a primitive root of unity using unitary bundles,
which gives a geometric explanation of the hyperbolic formulas?
\item Are there results for other values of $q$, for instance, $q$
negative real, or $q = \exp(\pi i/ r)$ but $j_{ab} \notin
[0,(r-2)/2]$?  Numerical experiments show that in some of these cases
the $6j$ symbols are rapidly increasing.
\item Are there similar formulas for $6j$ symbols related to the
Kashaev-Reshetikhin invariants \cite{kash:inv}?
\item There are analogs of the $6j$ symbols for other groups.  For
example, for any four-tuple of dominant weights
$j_{12},j_{23},j_{34},j_{41}$ such that the tensor products
$V_{j_{ab}} \otimes V_{j_{bc}}$ are multiplicity-free, one has two
canonical bases for the space of invariants $V_{j_{12}} \otimes
V_{j_{23}} \otimes V_{j_{34}} \otimes V_{j_{14}}$ given by the two
ways of pairing.  Are there explicit asymptotic formulas for these
symbols?
\end{enumerate}

\appendix

\section{Maple code}

The following code (not very well optimized) was used to generate
Figure 3.

\begin{Verbatim}[fontsize=\small]

with(plots);  with(linalg);
# quantum integer,factorial
pi:=evalf(Pi); nq:=(r,n)->sin(pi*n/r)/sin(pi/r);
faq:=(r,n)->evalf(product(nq(r,ii),ii=1..n),30); 
# quantum triangle symbol, Racah sum, quantum sixj symbol
trq:=(r,a,b,c)->faq(r,(b+c-a)/2)*faq(r,(a+c-b)/2)*
   faq(r,(a+b-c)/2)/(faq(r,(a+b+c)/2+1));
syq:=(r,a1,a2,a3,a4,b1,b2,b3)->sum((-1)^z*faq(r,z+1)/
   (faq(r,b1-z)*faq(r,b2-z)*faq(r,b3-z)*faq(r,z-a1)*
   faq(r,z-a2)*faq(r,z-a3)*faq(r,z-a4)),
   z=max(a1,a2,a3,a4)..min(b1,b2,b3));
sjq:=(r,a,b,c,d,e,f)->`if`(min(a+b-c,a-b+c,-a+b+c,a+e-f,e-a+f,-e+a+f,
   d+c-e,d-c+e,-d+c+e,b+f-d,b-f+d,-b+f+d,2*r-4-a-b-c,2*r-a-e-f,2*r-d-c-e,
   2*r-b-d-f)>= 0,evalf((trq(r,a,b,c)*trq(r,a,e,f)*trq(r,d,c,e)*
   trq(r,b,d,f))^(1/2)*syq(r,(a+b+c)/2,(a+e+f)/2,(c+d+e)/2,(b+d+f)/2,
   (a+b+d+e)/2,(a+c+d+f)/2,(b+c+e+f)/2),30),0);
# length corresponding to a dominant weight (integer)
len:=j->pi*(j+1)/r;
# The Gram matrix, its determinant and inverse
G:=(i,j,k,l,m,n)->linalg[matrix](4,4,[1, cos(i),cos(k),cos(m),
                               cos(i),1,     cos(j),cos(n),
                               cos(k),cos(j),1,     cos(l),
                               cos(m),cos(n),cos(l),1]);
detG:=(i,j,k,l,m,n)->evalf(det(G(i,j,k,l,m,n)),10);
Ginv:=(i,j,k,l,m,n)->evalf(inverse(G(i,j,k,l,m,n)));
# Does the tetrahedron exist?
tetexist:=(r,a,b,c,d,e,f)->`if`(min(a+b-c+1,a-b+c+1,-a+b+c+1,a+e-f+1,
   e-a+f+1,-e+a+f+1,d+c-e+1,d-c+e+1,-d+c+e+1,b+f-d+1,b-f+d+1,-b+f+d+1,
   2*r-4-a-b-c-3,2*r-a-e-f-3,2*r-d-c-e-3,2*r-b-d-f-3,
   detG(len(a),len(b),len(c),len(d),len(e),len(f))) >= 0,1,0);
# The amplitude of the asymptotic formula
amppredict:=(r,i,j,k,l,m,n)->evalf(((r/pi)^3*1.5*pi*(abs(detG(
   len(i),len(j),len(k),len(l),len(m),len(n)))/(36))^(1/2))^(-1/2),10);
# The Dihedral angles 
preangle:=(i,j,k,l,m,n)->evalf(seq(seq(pi - 
   arccos(-Ginv(i,j,k,l,m,n)[a,b]/(Ginv(i,j,k,l,m,n)[a,a]*
   Ginv(i,j,k,l,m,n)[b,b])^(1/2)),b=a+1..4),a=1..4),10);
diangle:=(r,i,j,k,l,m,n)->evalf(preangle(len(i),len(j),len(k),len(l),
   len(m),len(n)));
# The "Dehn invariant" aka Regge action 
dehn:=(r,i,j,k,l,m,n)-> diangle(r,i,j,k,l,m,n)[1]*(l+1)/2+
   diangle(r,i,j,k,l,m,n)[2]*(n+1)/2+ diangle(r,i,j,k,l,m,n)[3]*(j+1)/2+
   diangle(r,i,j,k,l,m,n)[4]*(m+1)/2+ diangle(r,i,j,k,l,m,n)[5]*(k+1)/2+
   diangle(r,i,j,k,l,m,n)[6]*(i+1)/2;
# The asymptotic formula for the degenerate case
G3:=(i,j,k)->linalg[matrix](3,3,[1,cos(i),cos(j),
                                cos(i),1,cos(k),
                                cos(j),cos(k),1]);
areaq:=(r,i,j,k)->evalf(det(G3(len(i),len(j),len(k))))^(1/2);
tang:=(r,i,j,k,l,m,n)->evalf(r^(-4/3)*2^(2/3)*3^(-1/3)*pi^(4/3)*
   GAMMA(2/3)^(-1)*(areaq(r,i,j,k)*areaq(r,i,m,n)*areaq(r,l,j,n)*
   areaq(r,l,m,k))^(-1/6));

# NUMERICAL EXPERIMENT 
r:=200; i:=40; j:=48; k:=50; l:=52; m:=54; stepsize:=.1; 
u:=2*max(i,j,k,l,m);
#calculate the 6j symbols
symbols:= [ seq([2*y,sjq(r,i,j,k,l,m,2*y)],y=0..u/2) ];
#does the terahedron exist
exist:= [seq([y*stepsize,tetexist(r,i,j,k,l,m,y*stepsize)],
   y=1..(r/2)/stepsize)];
#the predictions for the tangent case, plus or minus
tapredict:= [ seq([y,`if`(tetexist(r,i,j,k,l,m,y)=1,
   tang(r,i,j,k,l,m,y),undefined)],y=0..u/2) ];
mtapredict:= [ seq([y,-`if`(tetexist(r,i,j,k,l,m,y)=1,
   tang(r,i,j,k,l,m,y),undefined)],y=(u/2)..u) ];
#the dihedral angles 
angpredict :=   [ seq([y*stepsize,diangle(r,i,j,k,l,m,
      y*stepsize)[2]],y=0..u/stepsize) ];
for y from 0 to u/stepsize do
	phqpredict[y+1] := `if`(tetexist(r,i,j,k,l,m,y*stepsize)=1,
               phqpredict[y] + angpredict[y+1][2]*stepsize/2,
               `if`(tetexist(r,i,j,k,l,m,(y+1)*stepsize)=1,
               dehn(r,i,j,k,l,m,(y+1)*stepsize),0))   od;
#the predicted amplitudes , plus or minus
predicta:=[ seq([y*stepsize, `if`(tetexist(r,i,j,k,l,m,y*stepsize)=1,
	amppredict(r,i,j,k,l,m,y*stepsize),undefined)],y=1..u/stepsize) ];
mpredicta:=[ seq([y*stepsize,`if`(tetexist(r,i,j,k,l,m,y*stepsize)=1,
	-amppredict(r,i,j,k,l,m,y*stepsize),undefined)],y=1..u/stepsize) ];
#the prediction for the non-degenerate case 
predict:= [seq([stepsize*y,predicta[y][2]*cos(pi*(1/4) +
	phqpredict[y])],y=1..u/stepsize)];
# Show all the plots
plot([symbols,predict,mtapredict,tapredict,predicta,mpredicta],
   style=[point,line,line,line,line,line],
#   color=[red,blue,green,green,yellow,yellow],
#   linestyle=[SOLID,SOLID,DOT,DOT,SOLID,SOLID], 
   symbol=circle,symbolsize=5);

\end{Verbatim}

\def\cprime{$'$} \def\cprime{$'$} \def\cprime{$'$}
  \def\polhk#1{\setbox0=\hbox{#1}{\ooalign{\hidewidth
  \lower1.5ex\hbox{`}\hidewidth\crcr\unhbox0}}}

\end{document}